\documentclass[12pt,a4paper]{amsart}

\usepackage{amsfonts,color}
\usepackage{amsthm}
\usepackage{amssymb}
\usepackage{mathtools}
\usepackage{amsmath}
\usepackage{amscd}
\usepackage[latin2]{inputenc}
\usepackage{t1enc}
\usepackage[raggedrightboxes]{ragged2e}
\usepackage[mathscr]{eucal}
\usepackage{indentfirst}
\usepackage{graphicx}
\usepackage{enumerate}
\usepackage{graphics}
\usepackage{comment}
\usepackage[margin= 2.2cm]{geometry}
\usepackage{epstopdf} 
\usepackage{ragged2e}
\usepackage{xcolor}
\RequirePackage[colorlinks=true]{hyperref}

\newtheorem{thm}{Theorem}[section]
\newtheorem{prop}[thm]{Proposition}

\newtheorem{cor}[thm]{Corollary}
\newtheorem{lemma}[thm]{Lemma}
\newtheorem{Def}{Definition}
\newtheorem{remark}{Remark}
\theoremstyle{plain}
\theoremstyle{definition}

\definecolor{mycol}{rgb}{0,0,1}
\numberwithin{equation}{section}

\makeatletter
\@namedef{subjclassname@2020}{%
\textup{2020} Mathematics Subject Classification}
\makeatother

\begin{document}

\title{Maass Lifts of Half-Integral Weight Eisenstein Series and Theta Powers}

\author{Ajit Bhand, Karam Deo Shankhadhar, Ranveer Kumar Singh}

\address[Ajit Bhand]{Department of Mathematics, Indian Institute of Science Education and Research Bhopal,
Bhopal Bypass Road,  Bhauri,
Bhopal 462 066,
Madhya Pradesh, India}
\email{abhand@iiserb.ac.in, ajit.bhand@gmail.com}

\address[Karam Deo Shankhadhar]{Department of Mathematics, Indian Institute of Science Education and Research Bhopal,
Bhopal Bypass Road,  Bhauri,
Bhopal 462 066,
Madhya Pradesh, India}
\email{karamdeo@iiserb.ac.in, karamdeo@gmail.com}

\address[Ranveer Kumar Singh]{NHETC, Department of Physics and Astronomy, Rutgers University,
126 Frelinghuysen Rd., Piscataway NJ08855, USA}
\email{ranveersfl@gmail.com}




\begin{abstract}
In this paper, we explicitly construct mock modular forms whose shadows are Eisenstein series of arbitrary integral and half-integral weight, level and character at the cusps $\infty$ and $0$. As an application, we give explicit construction of harmonic weak Maass forms which are Hecke eigenforms and  are the preimages of 
$\Theta^k, k\in \{3,5,7\}$ under the shadow operator, where $\Theta$ is the classical Jacobi theta function. 
\end{abstract}

\maketitle

\section{Introduction}
\label{sec 1}

\noindent The theory of mock modular forms can be traced back to Ramanujan's famous ``deathbed" letter to G. H. Hardy written in 1920 in which he described certain functions which have ``modular like" properties. These mysterious functions, called \emph{mock theta functions} by Ramanujan, were not completely understood until the seminal work of Zwegers \cite{Z} as well as Bruinier and Funke \cite{BF} which provided the proper mathematical framework to study such functions. We now know that Ramanujan's  mock theta functions are the holomorphic parts of certain nonholomorphic functions called \emph{harmonic weak Maass forms} (see Section \ref{sec 2} for the definition). 

A mock modular form is defined as the holomorphic part of a harmonic weak Maass form. Given $k \in \frac{1}{2}\mathbb{Z}$, there exists a differential operator $\xi_{2-k}$, called the shadow operator,  that maps the space of harmonic weak Maass forms of weight $(2-k)$ surjectively to  the space of weakly holomorphic modular forms of weight $k$, and the image of a harmonic weak Maass form under this operator is called its \emph{shadow}. Given a weakly homomorphic modular form $f$ of weight $k$, it is natural to ask which harmonic weak maass forms have $f$ as the shadow. We consider these weak hamrmonic Maass forms as ``Maass lifts" of $f$.

Lifts of certain weight-$3/2$ weakly holomorphic modular forms have been constructed in \cite{Duke}.  Rhoades and Waldherr \cite{RW} have constructed a lift of of $\Theta^3$, where $\Theta$ is the classical Jacobi theta function. In \cite[Theorem 6.15]{Ono1}, Maass lifts of the standard Eisenstein series of integral and half-integral weight with trivial character at the cusp $\infty$ have been given. In \cite{W}, Wagner has constructed a lift of the Cohen-Eisenstein series which is a Hecke eigenform.  More recently, lifts of integral weight Eisenstein series with arbitrary level and character corresponding to the cusp $\infty$ and $\Theta^{2k}, k \in \{1,2,3,4\}$ have been constructed by Herrero and Pippich \cite{HP}. 

In this paper, we take a step forward and construct the Maass lifts of {general} Eisenstein series of half-integral weight  corresponding to the cusps $\infty$ and $0$ (Theorem \ref{thm 1.2}). As an application of the result, we construct lifts of $\Theta^5$ and $\Theta^7$ (Theorem \ref{thm 1.3}). Additionally, to make the treatment of lifts of Eisenstein series complete, we also construct lifts of integral weight Eisenstein series corresponding to the 
cusps $\infty$ and $0$ (Theorem \ref{thm 1.1}). In Theorem \ref{thm 1.4}, using our methods and a few results 
of \cite{RW} on the Kloosterman zeta-function, we construct a lift of $\Theta^3$, thereby providing a different proof of Rhoades and Waldherr's result. For $\Theta^k$, $k \in \{3, 5, 7\}$, the constructed preimages are Hecke eigenforms. Note that these are the odd powers of $\Theta$ which can be written as a linear combination of half-integral weight Eisenstein series (see \cite[Chapter IV]{NK}).      

We call the lifts of Eisenstein series that we have constructed as \emph{mock Eisenstein series.} The importance of these mock Eisenstein series lies in the following fact: Let $H_k^{\#}(\Gamma,\nu)$ denote the subspace of the space of harmonic weak Maass forms of weight $k$ with multiplier $\nu$ with respect to a finite index subgroup $\Gamma$  of $SL(2, \mathbb{Z})$ containing forms with at most polynomial growth at the cusps of $\Gamma$. All the preimages of Eisenstein series and certain powers of Jacobi theta functions constructed in this paper lie in this subspace. Zagier's famous weight-$3/2$ nonholomorphic Eisenstein series which is the preimage of Jacobi theta function and whose holomorphic part is the generating function of Hurwitz class numbers also lies in this subspace \cite{DMZ}. 

 In Section 5, we give a construction of general mock Eisenstein series corresponding to various cusps of $\Gamma$ which constitute a spanning set for $H_k^{\#}(\Gamma,\nu)$.
 While it is known that the space of harmonic weak Maass forms which have at most polynomial growth at every cusp is finite-dimensional and is spanned by the pre-images of Eisenstein series, to our knowledge it has not been written down in the existing literature. For this reason, we also provide a proof of this result (Theorem \ref{thm 1.5}).
  
The paper is organised as follows: In Section \ref{sec 1.2}, we state the results of this paper. In Section \ref{sec 2}, we recall the basic definitions of harmonic weak Maass forms and the corresponding Hecke theory and also set up our notations. In Section \ref{sec 3}, we review holomorphic Eisenstein series in the most general setup i.e. for a general finite index subgroup, arbitrary multiplier system and corresponding to an arbitrary cusp. In Section \ref{sec 4}, we review nonholomorphic Eisenstein series and analytic continuation of their Fourier coefficients. By using this analytic continuation, we define Eisenstein series of weights $2$ and $3/2$. 
In Section \ref{sec 5}, we construct our mock Eisenstein series. 
Finally, in Sections \ref{sec 6}, \ref{sec 7} and \ref{sec 8}, we provide the proofs of our results.

\subsection{Statement of results}
\label{sec 1.2}
To state the main results of the paper, we first establish some notations.  Let $\mathbb{C}$ be the complex plane. For each $z \in \mathbb{C}$, denote the real and imaginary parts of $z$ by $\text{Re}(z)$ and $\text{Im}(z)$, respectively. We also define $i= \sqrt{-1}$. Let $\mathbb{H}:= \{z \in \mathbb{C}:\; \text{Im}(z)> 0\}$ be the upper half-plane. For $z \in \mathbb{H}$, we let $q:=e^{2\pi i z}$. Let $k\in\frac{1}{2}\mathbb{Z}$ and $N\in\mathbb{N}$ with $4|N$ when $k\in\frac{1}{2}+\mathbb{Z}$. The congruence subgroup $\Gamma_0(N)$ of level $N$ is defined as 
\[
\Gamma_0(N) = \left\{\begin{pmatrix}
a&b\\c&d\\
\end{pmatrix}\in \mathrm{SL}(2,\mathbb{Z}):c~\equiv ~0~(\text{mod} ~N)\right\}.
\]
For any Dirichlet character $\chi$ modulo $N$, denote by $\chi^0$ the unique primitive character inducing $\chi$ and $m_{\chi}$ be the conductor of $\chi$ and put $\ell_{\chi}=N/m_{\chi}.$ We denote by $\mathbf{1}_N$ the principal character $\bmod~N$, and the trivial character which takes the value 1 on every integer is denoted simply by $\mathbf{1}$, and will be clear from the context.  Define the twisted divisor functions 
\begin{equation*}
\begin{split}
&\sigma_{k-1}^{\chi}(n)=\sum_{0<c \mid n} c^{k-1} \sum_{0<d \mid \operatorname{gcd}\left(\ell_{\chi}, c\right)} d \mu\left({\ell_{\chi}}/{d}\right) \overline{\chi^{0}}\left({\ell_{\chi}}/{d}\right) \chi^{0}\left({c}/{d}\right),\\
&\widehat{\sigma}_{1-k}^{\overline{\chi}}(n)=\sum_{\substack{c=1\\\text{gcd}(c,N)=1}}\overline{\chi}(c)c^{1-k},
\end{split}
\label{eq 1.1}
\end{equation*}
where $\mu$ denotes the M\"{o}bius function.
The Gauss sum of $\chi^0$ is defined by
\[
\tau(\chi^0)=\sum_{m\bmod m_{\chi}}\chi^0(m)e^{2\pi i\frac{m}{m_{\chi}}}.
\]
We denote by $$L(s,\chi)=\sum_{n=1}^{\infty}\frac{\chi(n)}{n^s},~\text{Re}(s)>1,$$ the Dirichlet $L$-function, and by $\zeta(s)=L(s, 1)$, the Riemann zeta function. 

Let $E_k(\Gamma_0(N),\chi,z)$ and $F_k(\Gamma_0(N),\chi,z)$ be the Eisenstein series corresponding to the cusps $\infty$ and $0$ respectively 
(see Section \ref{sec 3} for the precise definitions). We have the following theorem. 
\begin{thm}
\label{thm 1.1}
Let $2\leq k\in\mathbb{Z}$ with the assumption that $\chi\neq \mathbf{1}_N$ when $k=2$. The series 
\[
\sum_{n=0}^{\infty}a_{2-k}^+(n,\chi)q^n\quad \text{and}\;\;\sum_{n=0}^{\infty}b_{2-k}^+(n,\chi)q^n,\]
where 
\[
\begin{split}
&a_{2-k}^+(n,\chi)=\frac{(-2i)^{2-k}\pi n^{1-k}}{N^{k}(k-1)}\tau(\overline{\chi}^0)\sigma_{k-1}^{\overline{\chi}}(n),~~~\ \ \ b_{2-k}^+(n,\chi)=\frac{(-2i)^{2-k}\pi}{N^{k/2}(k-1)}\widetilde{\sigma}_{k-1}^{\overline{\chi}}(n),~~~n\neq 0,\\&a_{2-k}^+(0,\chi)=\frac{(-2i)^{2-k}\pi}{N^{k-1}(k-1)} \times \begin{cases}\zeta(k-1)\prod\limits_{p|N}(1-p^{-1})&\text{if $\chi=\mathbf{1}_N$}\\0&\text{if $\chi\neq\mathbf{1}_N$},
\end{cases},\\&b_{2-k}^+(0,\chi)=\frac{(-2i)^{2-k}\pi}{N^{k/2}(k-1)}L(k-1,\overline{\chi})
\end{split}
\]
are mock modular forms of weight $2-k$ for $\Gamma_0(N)$ with character $\chi$ whose shadows are  $E_k(\Gamma_0(N),\overline{\chi},z)$ and $F_k(\Gamma_0(N),\overline{\chi},z)$ respectively. 
\label{thm 1.1}
\end{thm}
\begin{remark}
 In \cite{HP}, Herrero and Pippich have constructed lifts of slightly more general integral weight Eisenstein series involving two Dirichlet characters corresponding to the cusp $\infty$. With certain simplifications, our construction of the lift of $E_k(\Gamma_0(N),\overline{\chi},z)$ coincides with that in \cite{HP}. 
\end{remark}

For the next result, define a new character 
\begin{equation}\label{eq:omega}
\omega(d)=\left(\frac{N}{d}\right)\chi(d),
\end{equation}
wher $\left(\frac{N}{d}\right)$ is the  Kronecker-Jacobi symbol. For the definition of the Kroncker-Jacobi symbol, we refer to \cite[Definition 3.4.6]{CS}.
With $j=\frac{2k-3}{2}$ and $n$ a non-zero integer, define primitive characters $\omega_{n}$ and $\omega^{1}$ by
\begin{equation}
\omega_{n}(a)=\left(\frac{-1}{a}\right)^{j}\left(\frac{nN}{a}\right)\omega(a) \quad \text { and } \quad \omega^{1}(a)=\omega(a)^{2} \text { for gcd}(a, n N)=1.\quad 
\label{eq 1.2}
\end{equation}
Also put 
\[
\begin{split}
&\beta_{2-k}(n,\omega)=\sum_{\substack{(ab)^2|n\\\text{gcd}(a,N)=1\\\text{gcd}(b,N)=1}} \mu(a) \omega_{n}(a) \omega^{1}(b) a^{\frac{1-2k}{2}}\\
&C_{2-k}(n,\chi)=\sum_{N|M| N^{\infty}}\left[\sum_{m\bmod~M}\left(\frac{M}{m}\right) \omega(m) \varepsilon_{m}^{-2k} e^{2\pi i\frac{nm}{M}} \right] M^{-k}
\end{split}
\]
where $N|M| N^{\infty}$ means that the sum runs over all $M$ defined as follows: for each $l\geq 1$, write $Nl=Md$ where $\text{gcd}(M,d)=1$ and $M$ divides all sufficiently high powers of $N$. These two requirements uniquely fix $M$ for a given $l$.  
Also, $\varepsilon_m=\sqrt{\left(\frac{-1}{m}\right)}$, where $\sqrt{\cdot}$ denotes the principal branch of the square root. Note that $\epsilon_m$ is $1$ if $m\equiv 1\bmod 4$ and 
$i$ when $m\equiv -1\bmod 4$.
We now state the first main result of the paper.
\begin{thm}
\label{thm 1.2}
Let $\frac{3}{2}\leq k\in\frac{1}{2}+\mathbb{Z}$ with the assumption that $\chi^2\neq \mathbf{1}_N$ 
when $k=3/2.$ Then the series 
\[
\sum_{n=0}^{\infty}A^+_{2-k}(n,\chi)q^n~~~\text{and}~~~\sum_{n=0}^{\infty}B_{2-k}^+(n,\chi)q^n,
\]
where 
\[
\begin{split}
&B^+_{2-k}(n,\chi)=\frac{(2i)^{2-k}\pi}{N^{k/2}(k-1)}\frac{L(\frac{2k-1}{2},\omega_n)}{L(2k-1,\omega^1)}\beta_{2-k}(n,\omega),~~~n\neq 0,\\
&B^+_{2-k}(0,\chi)=\frac{(2i)^{2-k}\pi}{N^{k/2}(k-1)}\frac{L\left(2k-2, \omega^{1}\right)}{L\left(2k-1, \omega^{1}\right)}
\end{split}
\]
and 
\[
A^+_{2-k}(n,\chi)=(-iN)^{k/2}B_{2-k}^+(n,\chi)C_{2-k}(n,\overline{\chi}, k-1)
\]
are mock modular forms of weight $2-k$ for $\Gamma_0(N)$ with character $\chi$ whose shadows are  $E_k(\Gamma_0(N),\overline{\chi},z)$ and $F_k(\Gamma_0(N),\overline{\chi},z)$ respectively.
\end{thm}

To state the next main result, we need to fix some notations. Write $n \neq 0$ as $n=f^{2} d$ with $d$ squarefree and $f=2^{q} w$ with $w$ odd. For $\kappa \in \{5, 7\}$, let
\begin{equation*}
c_{n}\left(\frac{\kappa-1}{2}\right)=\left\{\begin{array}{ll}
\frac{2^{\frac{\kappa-1}{2}}-\psi_{-n}(2)}{2^{\frac{\kappa-1}{2}}-1}, \hspace{1cm}n \equiv 1,2(\bmod ~4)&   \\&\\
\frac{1-\psi_{-n}(2) 2^{\frac{1-\kappa}{2}}}{1+2^{\frac{1-\kappa}{2}}}-2^{\frac{3-\kappa}{2}+Q(2-\kappa)}\left(\sigma_{\kappa-2}(2^Q)-2^{\frac{\kappa-3}{2}}\psi_{-n}(2)\sigma_{\kappa-2}(2^{Q-1})\right), & n \equiv 0, 3(\bmod ~ 4),\end{array}\right.
\label{eq:cn}
\end{equation*}
where
\begin{equation}
Q:=\left\{\begin{array}{ll}
q & \text { if } d \equiv 3(\bmod 4) \\
q-1 & \text { if } d \equiv 1,2(\bmod 4)
\end{array}\right.
\label{eq:Q}
\end{equation}
and $\psi_{n}(\cdot):=\left(\frac{D}{.}\right)$ with $D$ the discriminant of $\mathbb{Q}(\sqrt{n})$. For a character $\chi$, define the function $T_{s}^{\chi}$ by
\begin{equation}
T_{s}^{\chi}(w):=\sum_{a \mid w} \mu(a) \chi(a) a^{s-1} \sigma_{2 s-1}\left({w}/{a}\right)
\label{eq 1.3}
\end{equation}
where $\sigma_{\ell}$ denotes the $\ell$th divisor sum. Define 
\begin{equation*}
Z_n=\begin{cases}
e^{\frac{3 \pi i}{4}} \frac{6}{\pi^{2}} \log (2) & \text { if } n=0 \\
e^{\frac{3 \pi i}{4}} \frac{6}{\pi^{2}} \log (2) \frac{T_{1}^{\psi_1}(w)}{w} & \text { if } -n \text { is a square } \\
e^{\frac{3 \pi i}{4}}\frac{6}{\pi^{2}} L\left(1, \psi_{-n}\right) \frac{T_{1}^{\psi_{-n}}(w)}{w} c_{n}(1) & \text { otherwise, }
\end{cases}
\end{equation*}
and for $\kappa=5,7$ define
\begin{equation*}
Z_n\left(\frac{\kappa-1}{2}\right)=\left\{\begin{array}{ll}
e^{\frac{3 \pi i}{4}} \frac{\zeta(\kappa-2)}{\zeta(\kappa-1)} \frac{1-2^{-(\kappa-2)}-2^{-(\kappa-1)/2}}{1-2^{-\kappa+1}} & \text { if } n=0 \\
e^{\frac{3 \pi i}{4}}\frac{1}{\zeta(\kappa-1)}\frac{T_{(\kappa-1)/2}(w)}{w^{\kappa-2}}\zeta\left(\frac{\kappa-1}{2}\right)\left(1-2^{\frac{3-\kappa}{2}}\right) & \text { if } -n \text { is a square } \\
e^{\frac{3 \pi i}{4}} \frac{L\left(\frac{\kappa-1}{2}, \psi_{-n}\right)}{\zeta(\kappa-1)}  \frac{T_{(\kappa-1)/2}^{\psi_{-n}}(w)}{w^{\kappa-2}} \cdot c_{n}\left(\frac{\kappa-1}{2}\right) & \text { otherwise}.
\end{array}\right.
\end{equation*}
Finally, let
\[
Z'_n=\begin{cases}e^{\frac{3 \pi i}{4}} \frac{L\left(\frac{7}{2}, \psi_{-n}\right)}{\zeta(7)} w^{-6} T_{7/2}^{\psi_{-n}}(w) \frac{1-\psi_{-n}(2) 2^{-7/2}}{1-2^{-7}}&n\neq 0\\\frac{\zeta(5)}{\zeta(6)}&n=0.
\end{cases}
\]
As an application of Theorem \ref{thm 1.2}, we extend the result of Rhoades and Waldherr\cite{RW} to construct Maass lifts of $\Theta^5$ and $\Theta^7$. This is recorded in the following result.
\begin{thm}
\label{thm 1.3}
For $\kappa\in\{5,7\}$, the functions 
\[
F_{\Theta^{\kappa}}^+(z)=\sum_{n=0}^{\infty}c_{F_{\Theta^{\kappa}}}^+(n)q^n,
\]
where 
\[
c_{F_{\Theta^{\kappa}}}^+(n)=\frac{(-2i)^{2-\kappa/2}\pi}{(\kappa/2-1)2^{\kappa/2}}.\begin{cases}
Z_n\left(\frac{\kappa-1}{2}\right)&\kappa=5\\ & \\\overline{Z_{-n}\left(\frac{\kappa-1}{2}\right)}-2\overline{Z'_{-n}(\frac{\kappa}{2})}&\kappa=7,
\end{cases}
\]
are mock modular forms of weight $2-\kappa/2$ for $\Gamma_0(4)$ with trivial character whose shadows are 
$\Theta^{\kappa}$, respectively. Moreover, the corresponding harmonic Maass form $F_{\Theta^{\kappa}}(z)$ is a Hecke eigenform with eigenvalues $1+p^{2-\kappa}$, for any odd prime $p$.
\end{thm}

In \cite{RW}, Rhoades and Waldherr construct a  harmonic weak Maass form $F_{\Theta}$ whose shadow is $\Theta^3$ by considering certain non-holomorphic Poincar\'{e} series.
In the next result, we modify the construction given in Theorem \ref{thm 1.3} and use analytic continuation of the Fourier coefficients of our mock Eisenstein series, along with some preparatory results of \cite{RW}, to obtain a lift of $\Theta^3$. Our proof provides a way to understand how $F_{\Theta}$ arises from Maass lifts of Eisenstein series.
\begin{thm}
\label{thm 1.4}
The function 
\[
F_{\Theta^{3}}^+(z)=\sum_{n=0}^{\infty}\overline{Z_{-n}}q^n
\]
is a mock modular form of weight $1/2$ for $\Gamma_0(4)$ with trivial character whose shadow is $\Theta^{3}.$ Moreover the corresponding harmonic Maass form $F_{\Theta^{3}}(z)$ is a Hecke eigenforms with eigenvalue $1+{p}^{-1}$, for any odd prime $p$.
\end{thm}
Let $\Gamma<\mathrm{SL}(2,\mathbb{Z})$ be a finite index subgroup and let $\nu$ be a multiplier system for $\Gamma.$ Let $\mathfrak{E}_{k}(\Gamma,\nu)$ denote the Eisenstein subspace of the space of holomorphic modular forms, and $H_k^{\#}(\Gamma,\nu)$ be as in the introduction. Let $\mathcal{E}_{2-k}(\Gamma,\nu,z,\mathfrak{a})$ be the mock Eisenstein series corresponding to the cusp $\mathfrak{a}$ of $\Gamma$ (see Section \ref{sec 5} for details).  
We have the following result.
\begin{thm}
\label{thm 1.5}
For $2< k\in\frac{1}{2}\mathbb{Z}$, the space $H_{2-k}^{\#}(\Gamma,\nu)$ is spanned by the mock Eisenstein series $\mathcal{E}_{2-k}(\Gamma,\nu,z,\mathfrak{a})$ as $\mathfrak{a}$ varies over the inequivalent cusps of $\Gamma$ and $H_{2-k}^{\#}(\Gamma,\nu)\cong \mathfrak{E}_{k}(\Gamma,\overline{\nu}).$ Moreover, if  
$3\leq k\in\mathbb{Z}$ and $\chi$ is a Dirichlet character modulo $N$ then we have 
\[
\text{dim}\ H^{\#}_{2-k}(\Gamma_0(N),\chi)=\sum_{\substack{C \mid N\\\text{gcd}(C, N / C) \mid N / m_{\overline{\chi}}}} \phi(\operatorname{gcd}(C, N / C)),
\]  
where $m_{\overline{\chi}}$ is the conductor of $\overline{\chi}.$
\end{thm}

\section{Harmonic weak Maass forms and Hecke theory}
\label{sec 2}
\noindent Let $\Gamma$ be a finite index subgroup of $\mathrm{SL}(2,\mathbb{Z})$ and $\nu:\Gamma\longrightarrow\mathbb{C}$ be a multiplier system for $\Gamma$. 
For $\gamma=\left(\begin{smallmatrix}a&b\\c&d\end{smallmatrix}\right)\in\mathrm{SL}(2,\mathbb{Z})$, let $j(\gamma,z):=cz+d.$ 
Let $k \in \frac{1}{2}\mathbb{Z}$.
For a function $f:\mathbb{H}\longrightarrow\mathbb{C},$ define the weight-$k$ slash operator
$
(f|_{k}\gamma)(z):=j(\gamma,z)^{-k}f(\gamma z),
$
where $\gamma z=(az+b)(cz+d)^{-1}.$ A holomorphic function $f:\mathbb{H}\longrightarrow\mathbb{C}$ is called {\em weakly holomorphic modular form} (resp. {\em holomorphic modular form}, resp. {\em cusp form}) of weight $k$ for 
$\Gamma$ and multiplier system $\nu$ if $f|_k\gamma=\nu(\gamma)f$ for every $\gamma\in\Gamma$ and $f$ is meromorphic (resp. holomorphic, resp. vanishes) at the cusps of $\Gamma$. We use the following notation for the corresponding $\mathbb{C}$-vector spaces.
\[
\begin{split}
&M_k^!(\Gamma,\nu) ~:~\text{space of weakly holomorphic modular forms}\\
&M_k(\Gamma,\nu) ~:~\text{space of holomorphic modular forms}\\
&S_k(\Gamma,\nu)~:~\text{space of cusp forms}.
\end{split}
\]

Define the weight-$k$ hyperbolic Laplacian by
\begin{equation*}
\Delta_k = -y^2\left(\frac{\partial^2}{\partial x^2}+\frac{\partial^2}{\partial y^2}\right)+iky\left(\frac{\partial}{\partial x}+i\frac{\partial}{\partial y}\right)=-4y^2\frac{\partial}{\partial z}\frac{\partial}{\partial \bar{z}}+2iky\frac{\partial}{\partial \bar{z}}, \ x = {\rm Re}(z), y={\rm Im}(z).
\end{equation*}
\begin{Def}
A smooth function $f:\mathbb{H}\rightarrow\mathbb{C}$ is called a \emph{harmonic weak Maass form} of weight $k$ and multiplier $\nu$ for the group $\Gamma$ if 
\begin{enumerate}
\item $f|_k\gamma=\nu(\gamma)f$ for every $\gamma\in\Gamma.$
\item $\Delta_k(f)=0$.
\item There exists a polynomial $P_f(z)\in\mathbb{C}[q^{-1}]$ such that $f(z)-P_f(z)=O(e^{-\epsilon y})$ as $y\rightarrow\infty$ for some $\epsilon>0$. Similar conditions hold at other cusps. The polynomial $P_f(z)$ is called the {\it principal part of $f$}. 
\end{enumerate}
If condition (3) in the above definition is replaced by $f(z)=O(e^{\epsilon y})$, then $f$ is said to be a harmonic weak Maass form of \emph{manageable growth}. We denote the space of harmonic weak Maass forms of weight $k$, multiplier $\nu$ for the group $\Gamma$ by $H_k(\Gamma,\nu)$, and the space of harmonic weak Maass forms of manageable growth by $H_k^{!}(\Gamma,\nu)$.  
\end{Def}

Let $\Gamma_\infty$ be the stabilizer of the cusp $\infty$ in $\mathrm{SL}(2,\mathbb{Z})$. 
Suppose that the group $\Gamma$ contains $\Gamma_\infty$ and
the multiplier $\nu$ is trivial on $\Gamma_\infty$.
Any $f\in H_k^{!}(\Gamma,\nu), k\neq 1,$ has a Fourier series expansion \cite[Lemma 4.3]{BF} of the form  
\begin{equation}
f(z)=f(x+iy)=\sum\limits_{n>> -\infty}c_f^{+}(n)q^n+ c_f^{-}(0)y^{1-k}+\sum\limits_{\substack{n<< \infty\\n\neq 0}}c_f^{-}(n)\Gamma(1-k,-4\pi ny)q^n,
\label{eq 2.1}
\end{equation}
where $\Gamma(s,z)$ is the incomplete gamma function defined by
\begin{equation*}
\Gamma(s,z)\coloneqq\int\limits_z^{\infty}e^{-t}t^s\frac{dt}{t}.
\end{equation*}
The notation $\sum\limits_{n>> -\infty}$ means $\sum\limits_{n=\alpha_f}^{\infty}$ for some $\alpha_f\in\mathbb{Z}$, and analogously we write $\sum\limits_{n<<\infty}$ for $\sum\limits_{n=-\infty}^{\beta_f}$ for some 
$\beta_f\in\mathbb{Z}$. 
We call 
\begin{equation*}
f^+(z):=\sum\limits_{n>> -\infty}c_f^{+}(n)q^n
\end{equation*}
the \emph{holomorphic part} of $f$, and 
\begin{equation*}
f^-(z):=c_f^{-}(0)y^{1-k}+\sum\limits_{\substack{n<< \infty\\n\neq 0}}c_f^{-}(n)\Gamma(1-k,-4\pi ny)q^n
\end{equation*}
 the \emph{nonholomorphic part} of $f$. Moreover, if $f\in H_k(\Gamma,\nu)$ then $c_f^-(0)=0$ and the sum in nonholomorphic part runs only over negative integers.

Define the shadow operator by
\begin{equation*}
\xi_k\coloneqq 2iy^k\overline{\frac{\partial}{\partial\bar{z}}}.
\end{equation*} 
It is related to the weight $k$-Laplacian operator as follows.
\begin{equation}
\Delta_k=-\xi_{2-k}\circ\xi_{k}.
\label{eq 2.2}
\end{equation}
The image of any harmonic weak Maass form under the shadow operator is a weakly holomorphic modular form \cite[Theorem 5.10]{Ono1}. In particular, for $k\neq 1$, we have
\[
\xi_{2-k}(f)\in\begin{cases}M_{k}^!(\Gamma,\overline{\nu})&\text{if $f\in H_{2-k}^{!}(\Gamma,\nu)$}\\
S_k(\Gamma,\overline{\nu}) &\text{if $f\in H_{2-k}(\Gamma,\nu).$}
\end{cases}
\]
Furthermore, given the Fourier series expansion of $f$ as in \eqref{eq 2.1}, we have 
\begin{equation}
\xi_{2-k}(f(z))=\xi_{2-k}(f^-(z))=(k-1)\overline{c_f^-(0)}-(4\pi)^{k-1}\sum\limits_{n>>-\infty}\overline{c_f^-(-n)}n^{k-1}q^n.
\label{eq 2.3}
\end{equation}
Moreover, the map $\xi_{2-k}$ is surjective and the kernel of this map is $M_{2-k}^!(\Gamma,\nu)$  
\cite{Ono1, BF}. 
The image of $f \in H_{2-k}^{!}(\Gamma,\nu)$ under the operator $\xi_{2-k}$ is called the 
\emph{shadow} of $f$. 

\begin{Def}
A \emph{mock modular form} of weight $2-k$ is the holomorphic part $f^+$ of a harmonic weak Maass form of weight $2-k$ for which $f^-$ is nontrivial. The weakly holomorphic modular form $\xi_{2-k}(f)$ is called the shadow of the mock modular form $f^+$. 
\end{Def}
Let $\chi$ be a Dirichlet character modulo $N$ satisfying 
\[
\chi(-1)=\begin{cases}
(-1)^{k}&\text{if $k\in\mathbb{Z}$}\\
1,&\text{if $k\in\frac{1}{2}+\mathbb{Z}$.}
\end{cases}
\]
We will usually consider the case $\Gamma=\Gamma_0(N)$ and $\nu=\Psi_{k,\chi}$ defined as follows. 
For any $M=\left(\begin{smallmatrix}a&b\\c&d\end{smallmatrix}\right)\in \Gamma_0(N)$, we have 
\begin{equation}
\Psi_{k, \chi}(M):=\left\{\begin{array}{ll}
\chi(d) & \text { if } k \in \mathbb{Z}\\
\chi(d)\left(\frac{c}{d}\right) \varepsilon_{d}^{-2 k} & \text { if } k \in \frac{1}{2} + \mathbb{Z}.
\end{array}\right.
\label{eq 2.4}
\end{equation}
We assume that $4 | N$ if $k \in \frac{1}{2}+\mathbb{Z}$. We simply use the notation 
$\Psi_k$ for $\Psi_{k,\chi}$ if $\chi$ is trivial. In this setting, one can study Hecke theory for harmonic weak Maass forms. Here we recall some basic facts about Hecke operators $T(m), m \in \mathbb{N}$, and refer the reader to \cite[Chapter 7]{Ono1} for more details. We have the following proposition for the action of the Hecke operators on the Fourier coefficients.

\begin{prop}\emph{(\cite[Proposition 7.1]{Ono1})}. 
Let $k \in \frac{1}{2} \mathbb{Z}$. Suppose that $f(z) \in H_{k}^{!}\left(\Gamma_{0}(N), \Psi_{k,\chi}\right)$ with the Fourier series expansion given by \eqref{eq 2.1}. Then the following statements hold.
\begin{enumerate}
\item For $m \in \mathbb{N}$, we have that $f \mid T(m) \in H_{k}^{!}\left(\Gamma_{0}(N), \Psi_{k,\chi}\right)$.
\item If $k \in \mathbb{Z}, p \nmid N$ is a prime, then, with $\epsilon \in\{\pm\} (n\neq 0$ for $\epsilon=-)$, 
we have
$$
c_{f \mid T(p)}^{\epsilon}(n)=c_{f}^{\epsilon}(p n)+\chi(p) p^{k-1} c_{f}^{\epsilon}\left(\frac{n}{p}\right).
$$
Moreover, if $n=0$ and $\epsilon =-$ then we have
$$
c_{f \mid T(p)}^{-}(0)=\left(p^{k-1}+\chi(p)\right) c_{f}^{-}(0).
$$
\item If $k \in \frac{1}{2} \mathbb{Z} \backslash \mathbb{Z}, p \nmid N$ is a prime, then, with 
$\epsilon \in\{\pm\}(n \neq 0$ for $\epsilon=-)$, we have
$$
c_{f \mid T\left(p^{2}\right)}^{\epsilon}(n)=c_{f}^{\epsilon}\left(p^{2} n\right)+\chi^{*}(p)\left(\frac{n}{p}\right) p^{k-\frac{3}{2}} c_{f}^{\epsilon}(n)+\chi^{*}\left(p^{2}\right) p^{2 k-2} c_{f}^{\epsilon}\left(\frac{n}{p^{2}}\right)
$$
where $\chi{*}(n):=\left(\frac{(-1)^{k-\frac{1}{2}}}{n}\right) \chi(n) .$ If $n=0$ and $\epsilon=-$ then we have
$$
c_{f \mid T\left(p^{2}\right)}(0)=\left(p^{-2+2 k}+\chi{*}\left(p^{2}\right)\right) c_{f}^{-}(0).
$$
\end{enumerate}
\end{prop}

If $k \in -\frac{1}{2}\mathbb{N}$ then we have the following commutation relation of Hecke operators with the shadow operator \cite[Eq. (7.5)]{Ono1}.
\begin{equation}
p^{d(1-k)} \xi_{k}(f \mid T(p^{d}))=\xi_{k}(f) \mid T(p^{d}),
\label{eq 2.5}
\end{equation}
where
$
p \nmid N \ {\text{is a prime and}}\ 
d:=\left\{\begin{array}{ll}
1 & \text { if } k \in \mathbb{Z} \\
2 & \text { if } k \in \frac{1}{2}+\mathbb{Z}.
\end{array}\right.
$ 

\section{Holomorphic Eisenstein series}
\label{sec 3}
\noindent In this section, we recall holomorphic Eisenstein series of integral and half-integral weights, and 
write down their explicit Fourier series expansions. We follow the notations 
of the classical textbook by Iwaniec \cite{I}. Let $\frac{5}{2}\leq k\in\frac{1}{2}\mathbb{Z}$, 
$\Gamma\subset\mathrm{SL}(2,\mathbb{Z})$ be a finite index subgroup and $\nu$ be a multiplier system for $\Gamma$ satisfying $\nu(-I)=e^{-\pi ik}$ if $-I (=(\begin{smallmatrix}-1&0\\0&-1\end{smallmatrix}))\in\Gamma$.
Let $\mathfrak{a}$ be a cusp
of $\Gamma$ and let $\sigma_{\mathfrak{a}}\in\mathrm{SL}(2,\mathbb{R})$ be such that $\sigma_{\mathfrak{a}}\infty=\mathfrak{a}$ and $\sigma_{\mathfrak{a}}^{-1}\Gamma_{\mathfrak{a}}\sigma_{\mathfrak{a}}$ contains 
the matrix $\left(\begin{smallmatrix}1&1\\0&1\end{smallmatrix}\right)$. Here $\Gamma_{\mathfrak{a}}=\Gamma\cap \sigma_{\mathfrak{a}}\Gamma_{\infty}\sigma_{\mathfrak{a}}^{-1}$, the stabiliser of the cusp $\mathfrak{a}$ 
in $\Gamma$. 
Such a $\sigma_{\mathfrak{a}}$ is called the \emph{scaling matrix} for the cusp $\mathfrak{a}$. Suppose that $\nu$ is trivial on $\Gamma_{\mathfrak{a}}$. Such cusps are called \emph{singular} with respect to $\nu$. For 
$\gamma_1=\left(\begin{smallmatrix}a&b\\c&d\end{smallmatrix}\right),\gamma_2\in \mathrm{SL}(2,\mathbb{Z}),$ put  
\begin{equation*}
j(\gamma_1,z)=cz+d,~~~w(\gamma_1,\gamma_2)=j(\gamma_1,\gamma_2z)^{k}j(\gamma_2,z)^{k}j(\gamma_1\gamma_2,z)^{-k}.
\label{eq 3.1'}
\end{equation*} 

We define the holomorphic Eisenstein series corresponding to the cusp $\mathfrak{a}$ by 
\begin{equation}
E_{k}(\Gamma,\nu,z,\mathfrak{a}) = \sum_{\gamma \in \Gamma_{\mathfrak{a}} \backslash \Gamma} \overline{\nu}(\gamma)\overline{w}(\sigma_{\mathfrak{a}}^{-1},\gamma) j\left(\sigma_{\mathfrak{a}}^{-1} \gamma, z\right)^{-k}.
\label{eq 3.1}
\end{equation}
One can check that the above sum does not depend on the choice of the representatives of 
$\Gamma_{\mathfrak{a}} \backslash \Gamma$. Moreover, 
$E_{k}(\Gamma,\nu,z,\mathfrak{a})\in M_k(\Gamma,\nu)$ \cite[Proposition 3.1]{I}.

The Fourier expansion of these Eisenstein series at any singular cusp $\mathfrak{b}$ with scaling matrix 
$\sigma_{\mathfrak{b}}$ is given by \cite[Eq. 3.15, Section 3.2]{I}.
\begin{equation}
\left(E_{k}|_k\sigma_{\mathfrak{b}}\right)(\Gamma,\nu,z,\mathfrak{a})=\delta_{\mathfrak{a}\mathfrak{b}}+\sum_{n=1}^{\infty} \eta_{\mathfrak{a}\mathfrak{b}}(n)q^n,
\label{eq 3.2}
\end{equation} 
 where 
\begin{equation}
\eta_{\mathfrak{a}\mathfrak{b}}(n):=\left(\frac{2 \pi}{i}\right)^{k} \frac{n^{k-1}}{\Gamma(k)} \sum_{c>0} c^{-k} S_{\mathfrak{a}\mathfrak{b}}(\nu,n,c)
\label{eq 3.3}
\end{equation}
and $S_{\mathfrak{a}\mathfrak{b}}(\nu,n,c)$ is the Kloosterman sum defined by 
\begin{equation}
S_{\mathfrak{a}\mathfrak{b}}(\nu,n,c)=\sum_{\gamma=\left(\begin{smallmatrix}
a &b\\
c &d
\end{smallmatrix}\right) \in \Gamma_{\infty} \backslash \sigma_{\mathfrak{a}}^{-1} \Gamma \sigma_{\mathfrak{b}} / \Gamma_{\infty}} \overline{\nu}_{\mathfrak{ab}}(\gamma) e^{2\pi i\frac{nd}{c}}
\label{eq 3.4}
\end{equation}
with
\begin{equation}
\nu_{\mathfrak{ab}}(\gamma):=\nu\left(\sigma_{\mathfrak{a}} \gamma \sigma_{\mathfrak{b}}^{-1}\right) w\left(\sigma_{\mathfrak{a}}^{-1}, \sigma_{\mathfrak{a}} \gamma \sigma_{\mathfrak{b}}^{-1}\right) w\left(\gamma \sigma_{\mathfrak{b}}^{-1}, \sigma_{\mathfrak{b}}\right).
\label{eq 3.5}
\end{equation}
Here $\delta_{\mathfrak{a}\mathfrak{b}}=1$ if $\mathfrak{a}\sim\mathfrak{b} (\mathfrak{a}$ is equivalent to 
$\mathfrak{b}$) for $\Gamma$ and 0 otherwise. 
The sums $S_{\mathfrak{a}\mathfrak{b}}$ can be simplified greatly in some special cases. 

For the remainder of the section, we will assume that $\Gamma=\Gamma_0(N)$, 
where $4|N$ if $k\in\mathbb{Z}+\frac{1}{2}$, and $\nu=\Psi_{k,\chi}$. We consider the cusps $0$ and $\infty$ of $\Gamma_0(N)$ which are singular for the multiplier $\Psi_{k,\chi}$. For simplicity, we denote the corresponding Eisenstein series $E_{k}(\Gamma_0(N),\Psi_{k,\chi},z,\infty)$ by $E_k(\Gamma_0(N),\chi,z)$ and $E_{k}
(\Gamma_0(N),\Psi_{k,\chi},z,0)$ by $F_k(\Gamma_0(N),\chi,z)$. We now further simplify the Fourier expansions at the cusp $\infty$ given by \eqref{eq 3.2} of these two Eisenstein series of integral and half integral weights.  

\subsection{Integral weight Eisenstein series}\label{subsec:intweighteis}
Let $3\leq k\in\mathbb{Z}.$ The scaling matrix for the cusp $\infty$ is the identity matrix $I$. Using the double coset decomposition \cite[Proposition 2.7]{I}, we have
\[
\Gamma_0(N)=\Gamma_{\infty}\bigsqcup\limits_{c>0}\bigsqcup\limits_{d\bmod c}\Gamma_{\infty}\left(\begin{smallmatrix}\star&\star\\c&d\end{smallmatrix}\right)\Gamma_{\infty},
\] 
where $\left(\begin{smallmatrix}\star&\star\\c&d\end{smallmatrix}\right)\in\Gamma_0(N).$ It is clear that 
\[
\Gamma_{\infty} \backslash \Gamma_0(N)/ \Gamma_{\infty}=\{I\}\sqcup\left\{ \gamma=\left(\begin{smallmatrix}\star&\star\\c&d\end{smallmatrix}\right)\in\Gamma_0(N)~|~c>0, d\bmod c,\text{gcd}(c,d)=1\right\}.
\] 
Thus we see that the sum defining $S_{\infty\infty}(\Psi_{k,\chi},n,c)$ is empty if $c$ is not a multiple of $N$. Now, for any $\gamma=\left(\begin{smallmatrix}\star&\star\\c&d\end{smallmatrix}\right)\in\Gamma_0(N)$ and the multiplier $\Psi_{k,\chi}$ we have 
\[
\nu_{\infty\infty}(\gamma)=\chi(d).
\]
Using the above double coset decomposition, we get 
\begin{equation}
S_{\infty\infty}(\Psi_{k,\chi},n,c)=\begin{cases}
\sideset{}{^{*}}\sum\limits_{m\bmod c}\overline{\chi}(m)e^{2\pi i\frac{nm}{c}}&\text{if $N|c$}\\0&\text{otherwise},
\end{cases}
\label{eq:Sinfinfintweight}
\end{equation}
where $\sideset{}{^{*}}\sum$ means that the sum is only over those residue classes which are invertible modulo $c$, that is, gcd$(c,d)=1.$ Thus using \eqref{eq 3.3} and \eqref{eq 3.4}, we get
\begin{equation*}
E_{k}(\Gamma_0(N),\chi,z)=1+ \sum_{n=1}^{\infty}a_k(n,\chi)q^n,
\end{equation*}
where 
\[
\begin{split}
a_k(n,\chi)&:=\left(\frac{2 \pi}{i}\right)^{k} \frac{n^{k-1}}{\Gamma(k)}\sum_{\substack{c>0\\N|c}}c^{-k}\sum_{m\bmod c}\overline{\chi}(m)e^{2\pi i\frac{nm}{c}}\\&=\left(\frac{2 \pi}{i}\right)^{k} \frac{n^{k-1}}{\Gamma(k)}\sum_{c=1}^{\infty}(Nc)^{-k}\sum_{m\bmod Nc}\overline{\chi}(m)e^{2\pi i\frac{nm}{cN}}.
\end{split}
\]
Using \cite[Eq. 7.2.46]{miyake}, we have 
\[
\begin{split}
a_k(n,\chi)=\left(\frac{2 \pi}{i}\right)^{k} \frac{n^{k-1}}{\Gamma(k)}\sum_{c|n}(Nc)^{-k}c\sum_{m\bmod N}\overline{\chi}(m)e^{2\pi i\frac{nm}{cN}}.
\end{split}
\]
Let $\chi^0$ be a primitive Dirichlet character inducing $\chi$ and let $m_{\chi}$ be the conductor of $\chi$. Also put $\ell_{\chi}=N/m_{\chi}.$ Then using \cite[Eq. 7.2.47]{miyake}, we have 
\[
\begin{split}
a_k(n,\chi)&=\left(\frac{2 \pi}{iN}\right)^{k} \frac{1}{\Gamma(k)}\tau(\overline{\chi}^0)\sum_{0<c|n}c^{k-1}\sum_{0<d \mid \operatorname{gcd}\left(\ell_{\overline{\chi}}, c\right)} d \mu\left(\frac{\ell_{\overline{\chi}}}{d}\right) \overline{\chi}^{0}\left(\frac{\ell_{\overline{\chi}}}{d}\right) \chi^{0}\left(\frac{c}{d}\right),
\end{split}
\]
where $\tau(\chi)$ denotes the Gauss sum of $\chi$. Now define
\[
\sigma_{k-1}^{\overline{\chi}}(n)=\sum_{0<c|n}c^{k-1}\sum_{0<d \mid \operatorname{gcd}\left(\ell_{\overline{\chi}}, c\right)} d \mu\left(\frac{\ell_{\overline{\chi}}}{d}\right) \overline{\chi}^{0}\left(\frac{\ell_{\overline{\chi}}}{d}\right) \chi^{0}\left(\frac{c}{d}\right).
\]
Then we have 
\begin{equation}
a_k(n,\chi)=\left(-\frac{2 \pi i}{N}\right)^{k} \frac{\tau(\overline{\chi}^0)}{(k-1)!}\sigma_{k-1}^{\overline{\chi}}(n).
\label{eq:akint}
\end{equation}
In the above computation, we used the fact that $\overline{\overline{\chi}^{0}}=\chi^0$. 
Next we evaluate the Kloosterman sum for the cusp $0$. A choice of scaling matrix for this cusp is $\sigma_0=\left(\begin{smallmatrix}0&-\frac{1}{\sqrt{N}}\\\sqrt{N}&0\end{smallmatrix}\right)$. Using the double coset decomposition \cite[Proposition 2.7]{I}, we have
\[
\sigma_0^{-1}\Gamma_0(N)=\Gamma_{\infty}\bigsqcup\limits_{c>0}\bigsqcup\limits_{d\bmod c}\Gamma_{\infty}\left(\begin{smallmatrix}\star&\star\\c&d\end{smallmatrix}\right)\Gamma_{\infty},
\] 
where $\left(\begin{smallmatrix}\star&\star\\c&d\end{smallmatrix}\right)\in\sigma_0^{-1}\Gamma_0(N).$ It is easy to see that 
\[
\sigma_0^{-1}\Gamma_0(N)=\left\{\left(\begin{smallmatrix}
a\sqrt{N}&\frac{b}{\sqrt{N}}\\c\sqrt{N}&d\sqrt{N}
\end{smallmatrix}\right)~|~\left(\begin{smallmatrix}
-c&-d\\aN&b
\end{smallmatrix}\right)\in\Gamma_0(N)\right\}.
\]
Thus we have 
\[
\begin{split}
\Gamma_{\infty} \backslash \sigma_0^{-1}\Gamma_0(N)/ \Gamma_{\infty}=\left\{ \gamma=\left(\begin{smallmatrix}\star&\star\\c\sqrt{N}&d\sqrt{N}\end{smallmatrix}\right)\in\sigma_0^{-1}\Gamma_0(N)~|~c>0, d\bmod c,~\text{gcd}(c,N)=\text{gcd}(c,d)=1\right\}\\\sqcup\{I\}.
\end{split}
\]
Thus, for any $\gamma=\left(\begin{smallmatrix}\star&\star\\c\sqrt{N}&d\sqrt{N}\end{smallmatrix}\right)\in\Gamma_{\infty} \backslash \sigma_0^{-1}\Gamma_0(N)/ \Gamma_{\infty}$, we have 
\[
\nu_{0\infty}(\gamma)=\nu(\sigma_0\gamma)=\nu\left(\left(\begin{smallmatrix}-c&-d\\\star&\star\end{smallmatrix}\right)\right)=\overline{\chi}(-c),
\] 
and 
\begin{equation}
S_{0\infty}(\Psi_{k,\chi},n,c)=\begin{cases}
\chi(-\ell)\sideset{}{^{*}}\sum\limits_{m\bmod \ell}e^{2\pi i\frac{nm}{\ell}}&\text{if $c=\sqrt{N}\ell,$ gcd$(\ell,N)=1$}\\0&\text{otherwise}.
\end{cases}
\label{eq:S0infintweit}
\end{equation}
Using \eqref{eq 3.4} and \eqref{eq 3.5} we have 
\begin{equation*}
F_{k}(\Gamma_0(N),\chi,z)=\sum_{n=1}^{\infty}b_k(n,\chi)q^n,
\label{eq 3.8}
\end{equation*}
where 
\[
\begin{split}
b_k(n,\chi)&=\left(\frac{2 \pi}{i}\right)^{k} \frac{n^{k-1}}{\Gamma(k)}\sum_{\substack{c=1\\\text{gcd}(c,N)=1}}^{\infty}(c\sqrt{N})^{-k}\chi(-c)\sideset{}{^{*}}\sum_{m\bmod c}e^{2\pi i\frac{nm}{c}}\\&=\left(-\frac{2 \pi i}{\sqrt{N}}\right)^{k} \frac{n^{k-1}}{\Gamma(k)}\sum_{\substack{c=1\\\text{gcd}(c,N)=1}}^{\infty}c^{-k}\chi(-c)\sum_{\delta|\text{gcd}(c,n)}\mu\left(\frac{c}{\delta}\right)\delta.
\end{split}
\]
Interchanging the order of summations, we get 
\[
\begin{split}
b_k(n,\chi)&=\left(-\frac{2 \pi i}{\sqrt{N}}\right)^{k} \frac{n^{k-1}}{\Gamma(k)}\sum_{\delta|n}\delta\sum_{\substack{c=1\\\text{gcd}(c\delta,N)=1}}^{\infty}(c\delta)^{-k}\chi(-c\delta)\mu\left(\frac{c\delta}{\delta}\right)\\&=\left(-\frac{2 \pi i}{\sqrt{N}}\right)^{k} \frac{1}{\Gamma(k)}\sum_{\delta|n}\left(\frac{n}{\delta}\right)^{k-1}\chi(\delta)\sum_{\substack{c=1\\\text{gcd}(c\delta,N)=1}}^{\infty}c^{-k}\chi(-c)\mu\left(c\right)\\&=\left(-\frac{2 \pi i}{\sqrt{N}}\right)^{k} \frac{1}{\Gamma(k)}\sum_{\substack{\delta|n\\\text{gcd}(\delta,n)=1}}\delta^{k-1}\chi\left(\frac{n}{\delta}\right)\sum_{\substack{c=1\\\text{gcd}(c,N)=1}}^{\infty}c^{-k}\chi(-c)\mu\left(c\right).
\end{split}
\]
Now define 
\[
\widetilde{\sigma}_{k-1}^{\chi}(n)=\sum_{\substack{c|n\\\text{gcd}(c,n)=1}}\chi\left(\frac{n}{c}\right)c^{k-1}.
\]
Noting that $\chi(-1)=(-1)^k$ and 
\[
\frac{1}{L(\chi,s)}=\sum_{\substack{c=1\\\text{gcd}(c,N)=1}}^{\infty}c^{-s}\chi(c)\mu\left(c\right),\quad \text{Re}(s)\gg 0,
\]
we get 
\begin{equation}
b_k(n,\chi)=\left(\frac{2 \pi i}{\sqrt{N}}\right)^{k} \frac{L(\chi,k)^{-1}}{(k-1)!}\widetilde{\sigma}_{k-1}^{\chi}(n).
\label{eq:bkint}
\end{equation}

\subsection{Half-integral weight Eisenstein series}\label{subsect:half-integral}
We now consider $\frac{5}{2}\leq k\in\mathbb{Z}+\frac{1}{2}$. We will use the calculations of previous section with the new multiplier system. We first consider the cusp $\infty$ with scaling matrix $\sigma_{\infty}=I$. For any $I\neq \gamma=\left(\begin{smallmatrix}\star&\star\\c&d\end{smallmatrix}\right)\in \Gamma_{\infty}\backslash\Gamma_0(N)/\Gamma_{\infty}$, we have 
\[
\nu_{\infty\infty}(\gamma)=\chi(d)\left(\frac{c}{d}\right) \varepsilon_{d}^{-2 k},
\]  
and 
\[
 w\left(I,\gamma\right) w\left(\gamma,I\right)=1.
\]
Using the double coset decomposition, we have  
\[
S_{\infty\infty}(\Psi_{k,\chi},n,c)=\begin{cases}
\sideset{}{^{*}}\sum\limits_{m\bmod c}\overline{\chi}(m)\left(\frac{c}{m}\right) \varepsilon_{m}^{2 k}e^{2\pi i\frac{nm}{c}}&\text{if $N|c$}\\0&\text{otherwise},
\end{cases}
\]
where we used the fact that $\overline{\varepsilon_d^{-1}}=\varepsilon_d.$ Thus using \eqref{eq 3.4} and \eqref{eq 3.5}, we have  
\begin{equation*}
E_{k}(\Gamma_0(N),\chi,z)=1+ \sum_{n=1}^{\infty}A_k(n,\chi)q^n,
\end{equation*}
where 
\[
\begin{split}
A_k(n,\chi)&=\left(\frac{2\pi}{i}\right)^k\frac{n^{k-1}}{\Gamma(k)}\sum_{\substack{c=1\\N|c}}^{\infty}c^{-k}\sideset{}{^{*}}\sum_{m\bmod c}\overline{\chi}(m)\left(\frac{c}{m}\right) \varepsilon_{m}^{2 k}e^{2\pi i\frac{nm}{c}}\\&=\left(\frac{2\pi}{i}\right)^k\frac{n^{k-1}}{\Gamma(k)}\sum_{\substack{c=1\\N|c}}^{\infty}\Upsilon_n(c)c^{-k},
\end{split}
\]
where 
\[
\Upsilon_n(c):=\sum_{m\bmod ~c}\overline{\chi}(m)\left(\frac{c}{m}\right) \varepsilon_{m}^{2 k}e^{2\pi i\frac{nm}{c}}.
\]

We now consider the cusp $0$. To proceed, we need to evaluate $\nu_{0\infty}(\gamma)$ for $I\neq\gamma\in\Gamma_{\infty} \backslash \sigma_0^{-1}\Gamma_0(N)/ \Gamma_{\infty}.$ The following lemma will be useful.
\begin{lemma}
\begin{enumerate}[(i)]
\item For $z\in\mathbb{C}\setminus\mathbb{R}$ and $k\in\frac{1}{2}\mathbb{Z}$, we have
\[
(-z)^{k}=z^{k}(-i~\mathrm{sgn}(\mathrm{Im}(z)))^{2k}.
\] 
\item For any $z\in\mathbb{H}$ and $k\in\mathbb{Z}+\frac{1}{2}$, we have 
\[
\begin{split}
z^{-k}(az+b)^{k}=\begin{cases}
\left(a+\frac{b}{z}\right)^{k}&a>0\text{ or }b>0\\i^{4k}\left(a+\frac{b}{z}\right)^{k}&a<0\text{ and }b<0.
\end{cases}
\end{split}
\]
\item For any matrix $g=\begin{pmatrix}
a&b\\c&d
\end{pmatrix} \in \Gamma_0(4)$ with $a,b<0$ and $c,d>0$, and $\frac{5}{2}\leq k\in\frac{1}{2}\mathbb{Z}\setminus\mathbb{Z}$ we have that
\begin{equation*}
\left(\frac{-b}{a}\right)\varepsilon_a^{-2k}=\left(\frac{c}{d}\right)\varepsilon_d^{-2k}.
\end{equation*} 
\end{enumerate}
\label{lemma:negsqrootcomp}
\begin{proof}
%
Parts (i) and (ii) follow from direct computations. We prove (iii) by using the modularity of the Eisenstein series $\widetilde{F}_{k}\coloneqq (2z)^{-k}\widetilde{E}_{k}\left(-\frac{1}{4z}\right)$ for $\frac{5}{2}\leq k\in\frac{1}{2}\mathbb{Z}\setminus\mathbb{Z}$ where $\widetilde{E}_{k}$ is defined as: 
\begin{equation*}
\widetilde{E}_{k}(z)=\sum_{\gamma\in\Gamma_{\infty}\setminus\Gamma_0(4)}J(\gamma,z)^{-k},\quad J(\gamma,z)=\Psi_{k,\chi}(\gamma)j(\gamma,z),\quad\chi=\mathbf{1}.
\end{equation*} 
For a proof of modularity of $\widetilde{E}_{k}$ and $\widetilde{F}_{k}$, we refer the reader to $\S$2 of Chapter IV of \cite{NK}. We have $\widetilde{F}_{k}(gz)=J(g,z)^{k}\widetilde{F}_{k}(z)$. Using the definition of $\widetilde{E}_{k}$, we get
\begin{equation*}
\widetilde{F}_{k}(gz)=2^{-k}\left(\frac{az+b}{cz+d}\right)^{-k}\widetilde{E}_{k}\left(-\frac{1}{4\left(\frac{az+b}{cz+d}\right)}\right)=2^{-k}\left(\frac{az+b}{cz+d}\right)^{-k}\widetilde{E}_{k}\left(g'\left(-\frac{1}{4z}\right)\right),
\end{equation*} 
where $g'=\begin{pmatrix}
d&-c/4\\-4b&a
\end{pmatrix}$. Since we have 
$\widetilde{E}_{k}\left(g'\left(-\frac{1}{4z}\right)\right)=J\left(g',-\frac{1}{4z}\right)^{k}\widetilde{E}_{k}\left(-\frac{1}{4z}\right),$ we get 
\begin{equation*}
J(g,z)^kz^{-k}=J\left(g',-\frac{1}{4z}\right)^{k}\left(\frac{az+b}{cz+d}\right)^{-k}.
\end{equation*}
Substituting the expression for $J$ and using Lemma \ref{lemma:negsqrootcomp} (ii), we get
\begin{equation*}
\begin{split}
\left(\frac{c}{d}\right)\varepsilon_d^{-2k}(cz+d)^{k}z^{-k}&=\left(\frac{-4b}{a}\right)\varepsilon_a^{-2k}\left(\frac{b}{z}+a\right)^{k}\left(\frac{az+b}{cz+d}\right)^{-k}\\&=\left(\frac{-b}{a}\right)\varepsilon_a^{-2k}i^{4k}z^{-k}\left(az+b\right)^{k}\left(\frac{az+b}{cz+d}\right)^{-k},
\end{split}
\end{equation*}
Next we note that $(az+b)^{-1},(cz+d)\in \mathbb{H}$. Thus we have 
\[
\left(-\frac{1}{(-az-b)}\right)^{k}(cz+d)^k=i^{2k}\left(\frac{1}{(-az-b)}\right)^{k}(cz+d)^k=i^{2k}\left(-\frac{cz+d}{az+b}\right)^k=i^{4k}\left(\frac{cz+d}{az+b}\right)^k.
\]
Here we have used Lemma \ref{lemma:negsqrootcomp} (i) along with the fact that $(-az-b)^{-1}$ and $(cz+d)/(az+b)$ are in the lower half plane. 
This gives 
 \begin{equation*}
\begin{split}
\left(\frac{c}{d}\right)\varepsilon_d^{-2k}(cz+d)^{k}z^{-k}&=\left(\frac{-b}{a}\right)\varepsilon_a^{-2k}z^{-k}(cz+d)^{k}.
\end{split}
\end{equation*}
Since this holds for every $z\in\mathbb{H}$, we get our result.
\end{proof}
\end{lemma}

We know that the right hand side of \eqref{eq 3.1} is well defined in the sense that the summands in \eqref{eq 3.1} are independent of the choice of representatives. So we pick particular representatives which make our computations easier. To this end, note that a typical representative $I\neq\gamma\in\Gamma_{\infty} \backslash \sigma_0^{-1}\Gamma_0(N)/ \Gamma_{\infty}$ has the form 
\[
\gamma=\begin{pmatrix}
a\sqrt{N}&b/\sqrt{N}\\c\sqrt{N}&d\sqrt{N}
\end{pmatrix},\quad a,b,c,d\in\mathbb{Z},c>0, d\bmod c,\sigma_0\gamma\in\Gamma_0(N).
\]
Now for $\left(\begin{smallmatrix}1&n\\0&1\end{smallmatrix}\right)\in\Gamma_{\infty}$, we have
\[
\begin{pmatrix}1&n\\0&1\end{pmatrix}\begin{pmatrix}
a\sqrt{N}&b/\sqrt{N}\\c\sqrt{N}&d\sqrt{N}
\end{pmatrix}=\begin{pmatrix}
(a+cn)\sqrt{N}&(b+dnN)/\sqrt{N}\\c\sqrt{N}&d\sqrt{N}
\end{pmatrix}.
\]
Since $c,d>0$, thus we can choose $n$ large enough so that $a+cn$ and $b+dnN$ are positive. So we can choose the non-identity representatives $\gamma$ of $\Gamma_{\infty} \backslash \sigma_0^{-1}\Gamma_0(N)/ \Gamma_{\infty}$ of the form
\begin{equation*}
\gamma=\begin{pmatrix}
a\sqrt{N}&b/\sqrt{N}\\c\sqrt{N}&d\sqrt{N}
\end{pmatrix},\quad a,b,c,d\in\mathbb{Z}_{>0},d\bmod c,\sigma_0\gamma=\begin{pmatrix}
-c&-d\\aN&b
\end{pmatrix}\in\Gamma_0(N).
\label{eq:reps0posabcd}
\end{equation*}
With this choice of representatives, we have 
\[
\nu(\sigma_0\gamma)=\chi(d)\left(\frac{aN}{b}\right)\varepsilon_{b}^{-2k}=\overline{\chi}(-c)\left(\frac{d}{-c}\right)\varepsilon_{-c}^{-2k},
\] 
where we used Lemma \ref{lemma:negsqrootcomp} (iii) and the fact that $\Gamma_0(N)\subseteq\Gamma_0(4)$ in the half-integral weight case. Now using the fact that  $\left(\frac{d}{-c}\right)=\left(\frac{d}{c}\right)$ for $c,d>0$ and $\varepsilon_{-c}=i\varepsilon_{c}^{-1},$ we get
\[
\nu(\sigma_0\gamma)=i^{-2k}\overline{\chi}(-c)\left(\frac{d}{c}\right)\varepsilon_{c}^{2k}.
\]
Next, we compute the factor $w\left(\sigma_{0}^{-1}, \sigma_{0}\gamma\right) w\left(\gamma,I\right)$. First observe that $w(\gamma,I)=1.$ We have 
\[
w\left(\sigma_{0}^{-1}, \sigma_{0}\gamma\right)=j\left(\sigma_{0}^{-1}, \sigma_{0}\gamma z\right)^kj\left(\sigma_{0}\gamma, z\right)^kj\left(\gamma,z\right)^{-k}.
\]
Any $I\neq\gamma\in\Gamma_{\infty} \backslash \sigma_0^{-1}\Gamma_0(N)/ \Gamma_{\infty}$ can written as $\gamma=\left(\begin{smallmatrix}a&b\\c\sqrt{N}&d\sqrt{N}\end{smallmatrix}\right)$ with $c>0$ and $adN-bc=1$ so that $\sigma_0\gamma=\left(\begin{smallmatrix}-c&-d\\aN&b\end{smallmatrix}\right)$. Thus we have 
\[
\begin{split}
w\left(\sigma_{0}^{-1}, \sigma_{0}\gamma\right)&=\left[-\sqrt{N}\left(\frac{-cz-d}{aNz+b}\right)\right]^k\left(aNz+b\right)^k\left(c\sqrt{N}z+d\sqrt{N}\right)^{-k}\\&=\left(\frac{cz+d}{aNz+b}\right)^k\left(aNz+b\right)^k\left(cz+d\right)^{-k}\\&=\left[\frac{\frac{c}{aN}\left(z+\frac{b}{aN}\right)-\frac{bc}{(aN)^2}+\frac{d}{aN}}{z+\frac{b}{aN}}\right]^k\left(aNz+b\right)^k\left(cz+d\right)^{-k}\\&=\left[\frac{\frac{c}{aN}\left(z+\frac{b}{aN}\right)+\frac{1}{(aN)^2}}{z+\frac{b}{aN}}\right]^k\left(aNz+b\right)^k\left(cz+d\right)^{-k},
\end{split}
\]
where we used the fact that $(cz)^k=c^kz^k$ for any $c\geq 0,z\in\mathbb{C}$ and $k\in\frac{1}{2}\mathbb{Z}.$ Now using Lemma \ref{lemma:negsqrootcomp} (i) along with the fact that $a>0$, we get
\[
w\left(\sigma_{0}^{-1}, \sigma_{0}\gamma\right)=\left(aNz+b\right)^k\left(z+\frac{b}{aN}\right)^{-k}\left(\frac{1}{aN}(cz+d)\right)^{k}\left(cz+d\right)^{-k}=1.
\]
Noting that the multiplier $\Psi_{k,\chi}$ should satisfy $\Psi_{k,\chi}(-I)=e^{\pi ik}$, we have 
\[
\chi(-1)i^{-2k}=(-1)^{-k}\implies \chi(-1)=1.
\] 
Hence we finally have 
\[
S_{0\infty}(\Psi_{k,\chi},n,c)=\begin{cases}
i^{2k}\varepsilon_{\ell}^{-2 k}\chi(\ell)\sum\limits_{m\bmod \ell}\left(\frac{m}{\ell}\right) e^{2\pi i\frac{nm}{\ell}}&\text{if $c=\sqrt{N}\ell,$ gcd$(\ell,N)=1$}\\0&\text{otherwise}.
\end{cases}
\] 
This gives, for $\frac{5}{2}\leq k\in\mathbb{Z}+\frac{1}{2}$,
\begin{equation*}
F_{k}(\Gamma_0(N),\chi,z)=\sum_{n=1}^{\infty}B_k(n,\chi)q^n,
\label{eq 3.8}
\end{equation*}
where 
\begin{equation*}
B_k(n,\chi)=\left(-\frac{2 \pi i}{\sqrt{N}}\right)^{k} \frac{n^{k-1}}{\Gamma(k)}\sum_{\substack{c=1\\\text{gcd}(c,N)=1}}^{\infty}\widetilde{\Upsilon}_n(c)c^{-k},
\label{eq:coeffeishalfint0}
\end{equation*}
where
\begin{equation*}
\widetilde{\Upsilon}_n(c):=i^{2k}\chi(c)\varepsilon_{c}^{-2 k}\sum\limits_{m\bmod c}\left(\frac{m}{c}\right) e^{2\pi i\frac{nm}{c}}.
\label{eq:gausshalfint0}
\end{equation*}

\section{Nonholomorphic Eisenstein series} \label{sec 4}
\noindent In this section, we recall nonholomorphic Eisenstein series, write down their Fourier series expansions and then discuss analytic continuation of their Fourier coefficients. By using analytic continuation, we finally write 
down Eisenstein series of weights $2$ and $3/2$.
We use the notations of Section \ref{sec 3}. For a complex variable $s\in\mathbb{C}$ with $k+2\text{Re}(s)>2$, define the nonholomorphic Eisesntein series
\begin{equation*}
E_{k}(\Gamma,\nu,z,\mathfrak{a},s)=\sum_{\gamma \in \Gamma_{\mathfrak{a}} \backslash \Gamma} \overline{\nu}(\gamma)\overline{w}(\sigma_{\mathfrak{a}}^{-1},\gamma) j\left(\sigma_{\mathfrak{a}}^{-1} \gamma, z\right)^{-k}|j\left(\sigma_{\mathfrak{a}}^{-1} \gamma, z\right)|^{-2s}.
\label{eq 4.1}
\end{equation*}
It is easy to check using absolute and uniform convergence that 
\begin{equation*}
E_{k}(\Gamma,\nu,Mz,\mathfrak{a},s)=\nu(M)j(M,z)^{k}|j(M,z)|^{2s}E_{k}(\Gamma,\nu,z,\mathfrak{a},s),~~~~~\forall~~~M\in\Gamma.
\label{eq 4.2}
\end{equation*}
In particular, we see that $y^s E_{k}(\Gamma,\nu,Mz,\mathfrak{a},s)$ transforms like a modular form, where 
$y =$ Im$(z)$. We now calculate the Fourier expansion of $y^sE_{k}(\Gamma,\nu,Mz,\mathfrak{a},s)$ at any singular cusp $\mathfrak{b}$. To simplify notations, put 
\[
\pi(\gamma,z)=\overline{\nu}(\gamma)\overline{w}(\sigma_{\mathfrak{a}}^{-1},\gamma) j\left(\sigma_{\mathfrak{a}}^{-1} \gamma, z\right)^{-k}|j\left(\sigma_{\mathfrak{a}}^{-1} \gamma, z\right)|^{-2s}.
\]

We have
\[
\left((y^sE_{k})|_k\sigma_{\mathfrak{b}}\right)(\Gamma,\nu,z,\mathfrak{a})=j(\sigma_{\mathfrak{b}},z)^{-k}\text{Im}(\sigma_{\mathfrak{b}}z)^s\sum_{\gamma\in\Gamma_{\mathfrak{a}}\backslash\Gamma}\pi(\gamma,\sigma_{\mathfrak{b}}z).
\]
To proceed, we note that there is a bijection between the cosets $\Gamma_{\mathfrak{a}}\backslash\Gamma$ and $\Gamma_{\infty}\backslash\sigma_{\mathfrak{a}}^{-1}\Gamma\sigma_{\mathfrak{b}}$ given by the map
\[
\begin{split}
f:\Gamma_{\mathfrak{a}}\backslash\Gamma &\longrightarrow\Gamma_{\infty}\backslash\sigma_{\mathfrak{a}}^{-1}\Gamma\sigma_{\mathfrak{b}}\\
\gamma &\longmapsto \sigma_{\mathfrak{a}}^{-1}\gamma\sigma_{\mathfrak{b}}.
\end{split}
\]
With this observation, we get 
\[
\begin{split}
&\left((y^sE_{k})|_k\sigma_{\mathfrak{b}}\right)(\Gamma,\nu,z,\mathfrak{a})=j(\sigma_{\mathfrak{b}},z)^{-k}\text{Im}(\sigma_{\mathfrak{b}}z)^s\sum_{\gamma\in\Gamma_{\infty}\backslash\sigma_{\mathfrak{a}}^{-1}\Gamma\sigma_{\mathfrak{b}}}\pi(\sigma_{\mathfrak{a}}^{-1}\gamma\sigma_{\mathfrak{b}},\sigma_{\mathfrak{b}}z)\\&=j(\sigma_{\mathfrak{b}},z)^{-k}|j(\sigma_{\mathfrak{b}},z)|^{-2s}y^s\sum_{\gamma\in\Gamma_{\infty}\backslash\sigma_{\mathfrak{a}}^{-1}\Gamma\sigma_{\mathfrak{b}}}\overline{\nu}(\sigma_{\mathfrak{a}}^{-1}\gamma\sigma_{\mathfrak{b}})\overline{w}\left(\sigma_{\mathfrak{a}}^{-1}, \sigma_{\mathfrak{a}} \gamma{\sigma_{\mathfrak{b}}^{-1}}\right) j(\gamma\sigma_{\mathfrak{b}}^{-1}, \sigma_{\mathfrak{b}}z)^{-k}|j(\gamma\sigma_{\mathfrak{b}}^{-1}, \sigma_{\mathfrak{b}}z)|^{-2s}\\&=y^s\sum_{\gamma\in\Gamma_{\infty}\backslash\sigma_{\mathfrak{a}}^{-1}\Gamma\sigma_{\mathfrak{b}}}\overline{\nu}(\sigma_{\mathfrak{a}}^{-1}\gamma\sigma_{\mathfrak{b}})\overline{w}\left(\sigma_{\mathfrak{a}}^{-1}, \sigma_{\mathfrak{a}} \gamma{\sigma_{\mathfrak{b}}^{-1}}\right)\overline{w}\left(\gamma\sigma_{\mathfrak{b}}^{-1}, \sigma_{\mathfrak{b}}\right)j(\gamma,z)^{-k}|j(\gamma,z)|^{-2s}\\&=y^s\sum_{\gamma\in\Gamma_{\infty}\backslash\sigma_{\mathfrak{a}}^{-1}\Gamma\sigma_{\mathfrak{b}}}\overline{\nu}_{\mathfrak{a}\mathfrak{b}}(\gamma)j(\gamma,z)^{-k}|j(\gamma,z)|^{-2s},
\end{split}
\]
where we used the definition of $w$ and the fact that $w$ is $\pm 1$. By \cite[Proposition 2.7]{I}, we have
\[
\gamma\in\Gamma_{\infty}\backslash\sigma_{\mathfrak{a}}^{-1}\Gamma\sigma_{\mathfrak{b}}=\delta_{\mathfrak{a}\mathfrak{b}} \Gamma_{\infty} \bigsqcup_{c>0}\bigsqcup_{d(\bmod c)} \Gamma_{\infty}\left(\begin{smallmatrix}
\star & \star \\
c & d
\end{smallmatrix}\right) \Gamma_{\infty}
\]
so that we have 
\[
\begin{split}
\left((y^sE_{k})|_k\sigma_{\mathfrak{b}}\right)(\Gamma,\nu,z,\mathfrak{a})&=\delta_{\mathfrak{a}\mathfrak{b}}y^s+y^s\sum_{I\neq\gamma\in \Gamma_{\infty} \backslash \sigma_{\mathfrak{a}}^{-1} \Gamma \sigma_{\mathfrak{b}}}\overline{\nu}_{\mathfrak{a}\mathfrak{b}}(\gamma)j(\gamma,z)^{-k}|j(\gamma,z)|^{-2s}\\&=\delta_{\mathfrak{a}\mathfrak{b}}y^s+y^s\sum_{I\neq\gamma\in \Gamma_{\infty} \backslash \sigma_{\mathfrak{a}}^{-1} \Gamma \sigma_{\mathfrak{b}} / \Gamma_{\infty}}\overline{\nu}_{\mathfrak{ab}}(\gamma)I_{\gamma}(z),
\end{split}
\] 
where 
\[
I_{\gamma}(z)=\sum_{\gamma'\in\Gamma_{\infty}}j(\gamma\gamma',z)^{-k}|j(\gamma\gamma',z)|^{-2s}.
\]
Thus for any $\gamma=\left(\begin{smallmatrix}\star&\star\\c&d\end{smallmatrix}\right)$, we have 
\[
\begin{split}
I_{\gamma}(z)&=\sum_{n\in\mathbb{Z}}(c(z+n)+d)^{-k}|c(z+n)+d|^{-2s}\\&=\frac{1}{c^{k+2s}}\sum_{n\in\mathbb{Z}}\left(z+\frac{d}{c}+n\right)^{-k}\left|z+\frac{d}{c}+n\right|^{-2s}\\&=\frac{1}{c^{k+2s}}\sum_{n\in\mathbb{Z}}\left(\tau+n\right)^{-k}\left|\tau+n\right|^{-2s},
\end{split}
\]
where we have defined $\tau=z+d/c$. This step is justified as $c>0$ for any $\gamma=\left(\begin{smallmatrix}\star&\star\\c&d\end{smallmatrix}\right)\in\Gamma_{\infty} \backslash \sigma_{\mathfrak{a}}^{-1} \Gamma \sigma_{\mathfrak{b}} / \Gamma_{\infty}.$ Noting that $\text{Im}(\tau)=\text{Im}(z)=y$, we write $\tau=x+iy$. We now use the following form of Poisson summation formula
\[
\sum_{n\in\mathbb{Z}}f(x+n)=\sum_{n\in\mathbb{Z}}\widehat{f}(n)e^{2\pi inx},
\] 
where $f$ is a sufficiently nice function and $\widehat{f}$ is the Fourier transform of $f$ defined by 
\[
\widehat{f}(u)=\int_{-\infty}^{\infty}dx~f(x)e^{-2\pi iux}.
\]
For fixed $y$, we use Poisson summation formula for 
\[
f(x)=\left(\tau+n\right)^{-k}\left|\tau+n\right|^{-2s}=\left(x+iy+n\right)^{-k}\left|x+iy+n\right|^{-2s}.
\]
We get 
\[
\begin{split}
\sum_{n\in\mathbb{Z}}\left(\tau+n\right)^{-k}\left|\tau+n\right|^{-2s}&=\sum_{n\in\mathbb{Z}}\left(\int_{-\infty}^{\infty}dx~\left(x+iy\right)^{-k}\left|x+iy\right|^{-2s}e^{-2\pi inx}\right)e^{2\pi inx}\\&=\sum_{n\in\mathbb{Z}}\left(\int_{iy-\infty}^{iy+\infty}d\omega~\omega^{-k}|\omega|^{-2s}e^{-2\pi in\omega}\right)e^{2\pi in(iy)}e^{2\pi inx}\\&=\sum_{n\in\mathbb{Z}}\left(\int_{iy-\infty}^{iy+\infty}d\omega~\omega^{-k}|\omega|^{-2s}e^{-2\pi in\omega}\right)e^{2\pi in\tau}\\&=\sum_{n\in\mathbb{Z}}\left(\int_{iy-\infty}^{iy+\infty}d\omega~\omega^{-k}|\omega|^{-2s}e^{-2\pi in\omega}\right)e^{2\pi i\frac{nd}{c}}q^n,
\end{split}
\]
where we have made a change of integration variable to $x\to \omega=x+iy$ in the second step and substituted $\tau=z+d/c$ in third step. Thus we have 
\[
I_{\gamma}(z)=\sum_{n\in\mathbb{Z}}\left(\frac{1}{c^{k+2s}}e^{2\pi i\frac{nd}{c}}\right)\varrho^k_n(s,y)q^n,
\]
where 
\begin{equation*}
\begin{split}
\varrho_n^{k}(s,y)&=\int\limits_{iy-\infty}^{iy+\infty}d\omega~\omega^{-k}|\omega|^{-2s}e^{-2\pi in\omega}\\&=y^{-k-2s}e^{2\pi ny}\int\limits_{-\infty}^{\infty}(t+i)^{-k}(t^2+1)^{-s}e^{-2\pi inty}dt\\&=y^{-k-2s}e^{2\pi ny}\int\limits_{-\infty}^{\infty}(t+i)^{-k-s}(t-i)^{-s}e^{-2\pi inty}dt.
\end{split}
\end{equation*}
Here, the second expression has been obtained by substituting $\omega=(t+i)y$. By  \cite[Theorem 7.2.5]{miyake} we have
\begin{equation}
\varrho_n^{k}(s,y)=\begin{cases}
(-2\pi i)^{k}\pi^{s}\Gamma(k+s)^{-1}n^{k-1+s}y^{-s}\Omega(4\pi ny,k+s,s)&\text{if $n>0$}\\
i^{-k}(2\pi)\Gamma(k+s)^{-1}\Gamma(s)^{-1}\Gamma(k-1+2s)(2y)^{-k+1-2s}&\text{if $n=0$}\\
(2i)^{-k}\pi^s\Gamma(s)^{-1}(-n)^{s-1}y^{-k-s}e^{4\pi ny}\Omega(-4\pi ny,s,k+s)&\text{if $n<0$}\\
\end{cases}
\label{eq 4.4}
\end{equation}
where 
\begin{equation*}
\Omega(y,\alpha,\beta)\coloneqq \frac{y^{\beta}}{\Gamma(\beta)}\int_{0}^{\infty}e^{-yu}(u+1)^{\alpha-1}u^{\beta-1}du
\end{equation*}
with $y > 0$ and $\alpha,\beta\in\mathbb{C}$ with Re$(\beta) > 0$.
This function satisfies \cite[Lemma 7.2.4]{miyake}
\begin{equation}
\begin{split}
\Omega(z,1-\beta,1-\alpha)=\Omega(z,\alpha,\beta),\\\Omega(z,\alpha,0)=1.
\end{split}
\label{eq 4.6}
\end{equation}
Putting everything together, we get 
\begin{equation}
\left(y^sE_{k}|_k\sigma_{\mathfrak{b}}\right)(\Gamma,\nu,z,\mathfrak{a})=\delta_{\mathfrak{a}\mathfrak{b}}y^s+\sum_{n=-\infty}^{\infty}\widehat{\eta}_{\mathfrak{a}\mathfrak{b}}(n,s,y)q^n,
\label{eq 4.7}
\end{equation}
where 
\begin{equation}
\widehat{\eta}_{\mathfrak{a}\mathfrak{b}}(n,s,y)=y^s\varrho^k_n(s,y)\sum_{c>0}\frac{S_{\mathfrak{a}\mathfrak{b}}(\nu,n,c)}{c^{k+2s}},
\label{eq:hatetas}
\end{equation}
where $S_{\mathfrak{a}\mathfrak{b}}$ is as in \eqref{eq 3.4}. 

We recover the holomorphic Eisenstein series $E_{k}(\Gamma,\nu,z,\mathfrak{a})$ defined by \eqref{eq 3.1} considered in the previous section by evaluating 
$E_{k}(\Gamma,\nu,z,\mathfrak{a},s)$ at $s=0$.
\begin{prop}
For $\frac{5}{2}\leq k\in\frac{1}{2}\mathbb{Z},$ the nonholomorphic Eisenstein series $E_{k}(\Gamma,\nu,z,\mathfrak{a},s)$ is analytic at $s=0$ and $E_{k}(\Gamma,\nu,z,\mathfrak{a},0)=E_{k}(\Gamma,\nu,z,\mathfrak{a})$.
\label{prop:nonholholeisen}
\end{prop}

Now, we define weight-$2$ Eisenstein series for $\Gamma_0(N)$ with character $\chi$ corresponding to the cusps $\infty$ and $0$ by evaluating the analytic continuation of the corresponding nonholomorphic 
Eisenstein series of weight $2$ at $s=0$. 
To do this, we need to evaluate the coefficients $\widehat{\eta}_{\mathfrak{a}\infty}(n,s,y)$ for the cusps 
$\mathfrak{a}=\infty, 0$ at $s=0$ for all $n \in \mathbb{Z}$. By proceeding similar to the calculations of Section \ref{subsec:intweighteis}, we see that the Kloostermann zeta function
\[
\sum_{c>0}\frac{S_{\mathfrak{a}\infty}(\Psi_2,n,c)}{c^{2+2s}}
\]
is finite at $s=0$ for any $n\neq 0, \mathfrak{a}=\infty, 0$ (see \eqref{eq:akint} and \eqref{eq:bkint} and the subsequent calculations). 
Using the fact that $\varrho^2_{n}(0,y)=0$ for $n \leq 0$, we see that 
$\widehat{\eta}_{\mathfrak{a}\infty}(n,0,y)=0$, for any $n< 0$, $\mathfrak{a}=\infty,0$. 
For $n \geq 0$, the value of $\widehat{\eta}_{\mathfrak{a}\infty}(n,0,y)$ depends on whether $\chi$ is trivial or not. From the calculations in Proof of Theorem \ref{thm 1.1} (see \eqref{eq:mockak} and \eqref{eq:mockbk}), 
we find that $\widehat{\eta}_{\mathfrak{a}\infty}(0,s,y)$ has a pole at $s=0$ for $k=2$ when the character is trivial. This pole is cancelled by the zero of  $\Gamma(s)^{-1}$ coming from $\varrho^2_{0}(s,y)$ but gives a nontrivial nonholomorphic part in the Eisenstein series. Thus to get holomorphic Eisenstein series of weight $2$, we restrict ourselves to $\chi\neq \mathbf{1}_N$.   
\begin{Def}\label{def:weight2}
\normalfont
For $k=2$ and $\chi\neq \mathbf{1}_N$, define the weight-$2$ Eisenstein series  $E_{2}(\Gamma_0(N),\chi,z)$ and $F_{2}(\Gamma_0(N),\chi,z)$ corresponding to the cusps $\infty$  and $0$ as follows. 
\[
E_{2}(\Gamma_0(N),\chi,z)=1+\sum_{n=1}^{\infty}a_2(n,\chi)q^n,\quad F_{2}(\Gamma_0(N),\chi,z)=\sum_{n=1}^{\infty}b_2(n,\chi)q^n,
\]
where $a_2(n,\chi)$ and $b_2(n,\chi)$ are given by \eqref{eq:akint} and \eqref{eq:bkint} respectively evaluated 
at $k=2$.
\end{Def} 
\begin{remark}
We see that the nonholomorphic Eisenstein series $E_{2}(\mathrm{SL}(2,\mathbb{Z}),\mathbf{1},z)$ coincides with $E_2^*(z):=E_2(z)-3/\pi y$, where $E_2(z)$ is the holomorphic weight-$2$ Eisenstein series for 
$\mathrm{SL}(2,\mathbb{Z})$. See \cite[Proof of Lemma 6.2]{Ono1} for more details.
\end{remark}


\subsection{Analytic continuation of Fourier coefficients}\label{subsec:anacontcoefhalfint}
\label{subsec 4.2}

Let $k\in\frac{1}{2}+\mathbb{Z}.$ In this subsection, we recall the analytic continuation of the Fourier coefficients of $E_{k}(\Gamma_0(N),\Psi_{k,\chi},z,\infty,s)$ and $E_{k}(\Gamma_0(N),\Psi_{k,\chi},z,0,s)$ at the cusp 
$\infty$. Proceeding similar to the calculations done in Section \ref{subsect:half-integral}, for 
$k+2 {\rm Re}(s) > 2$ we have 
\begin{equation}
\begin{split}
&E_{k}(\Gamma_0(N),\Psi_{k,\chi},z,\infty,s)=1+\sum_{n\in\mathbb{Z}}A_{k}(n,\chi,s)\varrho^k_n(s,y)q^n,\\
&E_{k}(\Gamma_0(N),\Psi_{k,\chi},z,0,s)=i^{2k}N^{-k/2-s}\sum_{n\in\mathbb{Z}}B_{k}(n,\chi,s)\varrho^k_n(s,y)q^n,
\end{split}
\label{eq:eisnhEF}
\end{equation}
where
\[
A_{k}(n,\chi,s):=\sum_{\substack{c=1\\N|c}}^{\infty}\Upsilon_n(c)c^{-k-2s},~~~
B_{k}(n,\chi,s):=i^{-2k}\sum_{\substack{c=1\\\text{gcd}(c,N)=1}}^{\infty}\widetilde{\Upsilon}_n(c)c^{-k-2s}.
\]
Writing out more explicitly, we have 
\begin{equation}
\begin{split}
& A_{k}(n,\chi,s)=\sum_{\substack{c=1\\N|c}}^{\infty}\left(\sum_{m (\bmod ~c)}\overline{\chi}(m)\left(\frac{c}{m}\right)\varepsilon_{m}^{2k}e^{2\pi i\frac{nm}{c}}\right)c^{-k-2s},\\& B_{k}(n,\chi,s)=\sum_{\substack{c=1\\\text{gcd}(c,N)=1}}^{\infty}\left(\chi(c)\varepsilon_{c}^{-2k}\sum_{m (\bmod ~c)}\left(\frac{m}{c}\right)e^{2\pi i\frac{nm}{c}}\right)c^{-k-2s}.
\end{split}
\label{eq:coeffnonholeisboth}
\end{equation}

Let $\omega$ be the character modulo $N$ defined by \eqref{eq:omega}.
We have the following theorem originally due to Shimura and Sturm.
\begin{thm}\emph{(\cite[Theorem 2, Lemma 2]{St})}
\label{thm 4.1}
Let $j=\frac{1-2k}{2},n$ a non-zero integer, and let $\omega_{n}, \omega^{1}$ be the characters defined by
\eqref{eq 1.2}.
Then
\[
\begin{split}
L\left(4s-2 j, \omega^{1}\right) B_k(n,\chi,s)=L\left(2s-j, \omega_{n}\right) \beta(n, 2s, \omega) \\
\beta(n, 2s, \omega)=\sum \mu(a) \omega_{n}(a) \omega^{1}(b) a^{j-2s} b^{2-2k-2s},
\end{split}
\]
where the last sum is extended over all positive integers $a$ and $b$ prime to $N$ such that $(a b)^{2}$ divides $n.$ If $n=0$ we have
\begin{equation}\label{eq:Bkanacont}
B_k(0,\chi,s)=\frac{L\left(4s+2k-2, \omega^{1}\right)}{L\left(4s-2j, \omega^{1}\right)}.
\end{equation}
For each $l\geq 1$, write $Nl=Md$ where $\text{gcd}(M,d)=1$ and $M$ divides all sufficiently high powers of $N$.  Let $N|M|N^{\infty}$, signify this condition. Then we have 
\begin{equation}\label{eq:Akanacont}
A_k(n,\chi,s)=B_k\left(n,\overline{\chi},s\right) C_k(n,\overline{\chi}, s),
\end{equation} 
where
$$
C_k(n,\overline{\chi},s)=\sum_{N|M| N^{\infty}}\left[\sum_{m\bmod~M}\left(\frac{M}{m}\right) \overline{\chi}(m) \varepsilon_{m}^{2k} e^{2\pi i\frac{nm}{M}} \right] M^{-k-2s}
$$
and $C_k(n,\overline{\chi},s)$ is a finite Dirichlet series if $n \neq 0.$
\label{thm:anacontAB}
\end{thm}

For $\frac{5}{2} \leq k \in \mathbb{Z}+\frac{1}{2}$, the coefficients $A_{k}(n,\chi,s)$ and $B_{k}(n,\chi,s)$ given by \eqref{eq:coeffnonholeisboth} are well-defined at $s=0$ and in this case we obtain holomorphic Eisenstein series as in Proposition \ref{prop:nonholholeisen}. However, for $k=3/2$, we have to use Theorem \ref{thm:anacontAB} and evaluate the Fourier coefficients  $A_{k}(n,\chi,s)$ and $B_{k}(n,\chi,s)$ at $s=0$. In particular, we need to take care of the possible poles of $L(1,\omega_n)$ and $L(1,\omega^1)$. It is clear that $L(1,\omega_n)$ is finite for every $n$ under the assumption that $\chi$ is not the inverse of the character defined by $\left(\frac{-n}{\cdot}\right)$, and $L(1,\omega^1)$ is finite under the assumption that
$\chi^2\neq \mathbf{1}_N$. 
Therefore we make the following definition.
\begin{Def}\label{def:eis3/2anacont}
\normalfont
Let $\chi$ be a Dirichlet character modulo $N$ such that $\chi^2\neq \mathbf{1}_N$. Define the Eisenstein series  
$E_{3/2}(\Gamma_0(N),\chi,z)$ and $F_{3/2}(\Gamma_0(N),\chi,z)$ of weight $3/2$ 
corresponding to the cusps $\infty$  and $0$ respectively as follows.
\[
\begin{split}
&E_{3/2}(\Gamma_0(N),\chi,z)=E_{3/2}(\Gamma_0(N),\Psi_{3/2},z)=1+\frac{4\sqrt{2}\pi}{i^{3/2}}\sum_{n=1}^{\infty}A_{3/2}(n,\chi)n^{1/2}q^n\\
&F_{3/2}(\Gamma_0(N),\chi,z)=F_{3/2}(\Gamma_0(N),\Psi_{3/2},z)=\frac{4\sqrt{2}\pi i^{3/2}}{N^{3/4}}\sum_{n=1}^{\infty}B_{3/2}(n,\chi)n^{1/2}q^n,
\end{split}
\]
where $A_{3/2}(n,\chi)$ and $B_{3/2}(n,\chi)$ are defined by evaluating the analytic continuations of 
$A_{3/2}(n,\chi,s)$ and $B_{3/2}(n,\chi,s)$ respectively at $s=0$ in Theorem \ref{thm:anacontAB}.
\end{Def}

\begin{remark}
For $\chi^2=\mathbf{1}_N$, $L(1,\omega_n)$ and $L(1,\omega^1)$ may have poles. Thus the coefficient $B_{3/2}(n,\mathbf{1}_N,s)$, and hence $A_{3/2}(n,\mathbf{1}_N,s)$ as well, 
$($see \eqref{eq:Bkanacont} and \eqref{eq:Akanacont}$)$ are singular at $s=0$ for some $n$.
These poles are cancelled by the zeros of $\varrho_n^{3/2}(0,y)$, but as a result some of the coefficients of Eisenstein series end up having the incomplete gamma function coming from 
$\varrho_n^{3/2}(0,y)$. Thus we end up by getting nonholomorphic Eisenstein series having correct modular transformations. In Section \ref{sec 7}, we show that a certain combination of the nonholomorphic Eisenstein series $E_{3/2}(\Gamma_0(4),\mathbf{1}_4,z)$ and $F_{3/2}(\Gamma_0(4),\mathbf{1}_4,z)$ 
gives us $\Theta^3$. This fact is crucial in our approach towards constructing the Maass lift of $\Theta^3$ in Theorem \ref{thm 1.4}.  
\end{remark}

\section{Mock Eisenstein series}\label{sec 5}

\noindent For $k\in\frac{1}{2}\mathbb{Z}$ and $s\in\mathbb{C}\backslash\{0\}$ with $\text{Re}(s)>k/2$, we define the function $\mathcal{E}_{2-k}(\Gamma,\nu,z,\mathfrak{a},s)$ by 
\begin{equation}
\mathcal{E}_{2-k}(\Gamma,\nu,z,\mathfrak{a},s)=\frac{y^s}{s}E_{2-k}(\Gamma,\nu,z,\mathfrak{a},s), \ \ y = \text{Im}(z).
\label{eq:nonholmock}
\end{equation}
One can easily check that $\mathcal{E}_{2-k}(\Gamma,\nu,z,\mathfrak{a},s)$ satisfies the following transformation property.
\begin{equation}
\mathcal{E}_{2-k}(\Gamma,\nu,\gamma z,\mathfrak{a},s)=\nu(\gamma)j(\gamma,z)^{2-k}\mathcal{E}_{2-k}(\Gamma,\nu,z,\mathfrak{a},s),~~~~\forall~~\gamma\in\Gamma.
\label{eq 4.10}
\end{equation}

For $k> 2$, the point $s=k-1$ is in the domain of the definition of $\mathcal{E}_{2-k}(\Gamma,\nu,z,\mathfrak{a},s)$. 
We define the {\it mock Eisenstein series} $\mathcal{E}_{2-k}(\Gamma,\nu,z,\mathfrak{a})$ 
of weight $2-k$ as follows. 
\begin{equation*}
\mathcal{E}_{2-k}(\Gamma,\nu,z,\mathfrak{a}):=\mathcal{E}_{2-k}(\Gamma,\nu, z,\mathfrak{a},k-1).
\label{eq 4.8}
\end{equation*}

\begin{thm}
\label{thm 4.2}
For $\frac{5}{2}\leq k\in\frac{1}{2}\mathbb{Z},$ we have that $\mathcal{E}_{2-k}(\Gamma,\nu, z,\mathfrak{a})\in H_{2-k}(\Gamma,\nu)$ with shadow $E_{k}(\Gamma,\overline{\nu},z,\mathfrak{a}).$ Moreover, $\mathcal{E}_{2-k}(\Gamma,\nu, z,\mathfrak{a})$ is bounded by a constant at every cusp of $\Gamma$ except at $\mathfrak{a}$ where it grows polynomially.

\begin{proof}
The modularity of $\mathcal{E}_{2-k}(\Gamma,\nu, z,\mathfrak{a})$ follows from \eqref{eq 4.10}. 
Consider the Fourier expansion of $\mathcal{E}_{2-k}(\Gamma,\nu,z,\mathfrak{a})$ at any singular cusp $\mathfrak{b}$.
\[
\left(\mathcal{E}_{2-k}|_{2-k}\sigma_{\mathfrak{b}}\right)(\Gamma,\nu,\gamma z,\mathfrak{a})=\frac{y^{k-1}}{k-1}\delta_{\mathfrak{a}\mathfrak{b}}+\frac{1}{k-1}\sum_{n=-\infty}^{\infty}\widehat{\eta}_{\mathfrak{a}\mathfrak{b}}(n,k-1,y)q^n.
\]
Next observe that for $n<0$, we have
\begin{equation*}
y^{\alpha-2}\Omega(-4\pi ny,\alpha-1,1)=-4\pi ny^{\alpha-1}\int_{0}^{\infty}e^{4\pi nyx}(x+1)^{\alpha-2}dx.
\end{equation*}
Substituting $t=x+1$ we get 
\begin{equation*}
y^{\alpha-2}\Omega(-4\pi ny,\alpha-1,1)=-4\pi ne^{-4\pi ny}y\int_{1}^{\infty}e^{4\pi nyt}(yt)^{\alpha-2}dt.
\end{equation*}
Again substituting $w=-4\pi nyt$ we get
\begin{equation}
y^{\alpha-2}\Omega(-4\pi ny,\alpha-1,1)=(-4\pi n)^{2-\alpha}e^{-4\pi ny}\Gamma(\alpha-1,-4\pi ny).
\label{eq:omeinga}
\end{equation}
Using this relation along with \eqref{eq 4.4} we have for $n<0$ 
\begin{equation}
\frac{y^{k-1}}{k-1}\varrho^{2-k}_n(k-1,y)=-\frac{(2\pi i)^{k}}{(4\pi)^{k-1}\Gamma(k)}\Gamma(k-1,-4\pi ny).
\label{eq:rho-}
\end{equation}
Next using \eqref{eq 4.6} we have for $n>0$
\begin{equation}
\frac{y^{k-1}}{k-1}\varrho^{2-k}_n(k-1,y)=\frac{1}{k-1}\left(-2i\right)^{2-k}\pi.
\label{eq:rho+}
\end{equation}
Finally 
\begin{equation}
\frac{y^{k-1}}{k-1}\varrho^{2-k}_0(k-1,y)=\frac{\pi}{k-1}(-2i)^{2-k}.
\label{eq:rho0}
\end{equation} 
Hence the Fourier expansion can be written in the cannonical way as 
\begin{equation}
\left(\mathcal{E}_{2-k}|_{2-k}\sigma_{\mathfrak{b}}\right)(\Gamma,\nu, z,\mathfrak{a})=\sum_{n=0}^{\infty}c_{2-k}^+(n)q^n+c_{2-k}^-(0)y^{k-1}+\sum_{n<0}c_{2-k}^-(n)\Gamma(k-1,-4\pi ny)q^n,
\label{eq 4.11}
\end{equation}
where 
\begin{equation}
c_{2-k}^+(n)=\frac{\pi}{k-1}(-2i)^{2-k}\sum\limits_{c>0}\frac{S_{\mathfrak{a}\mathfrak{b}}(\nu,n,c)}{c^{k}};~~~~~~~c_{2-k}^-(n)=\begin{cases}
-\frac{(2\pi i)^k}{(4\pi)^{k-1}\Gamma(k)}\sum\limits_{c>0}\frac{S_{\mathfrak{a}\mathfrak{b}}(\nu,n,c)}{c^{k}}&\text{if $n<0$}\\
\frac{1}{k-1}\delta_{\mathfrak{a}\mathfrak{b}}&\text{if $n=0$}.
\end{cases}
\label{eq 4.12}
\end{equation}
The above Fourier series expansion implies the claimed cusp conditions. 

Now we compute the shadow of $\mathcal{E}_{2-k}(\Gamma,\nu, z,\mathfrak{a})$. 
Note that if $\nu$ is a multiplier of weight $2-k$, then $\overline{\nu}$ is a multiplier of weight $k-2$ and hence also a multiplier of weight $k$ (cf. \cite[page 42]{I}). One can easily see by using \eqref{eq 2.3} that 
\[
\xi_{2-k}\left(\left(\mathcal{E}_{2-k}|_{2-k}\sigma_{\mathfrak{b}}\right)(\Gamma,\nu,z,\mathfrak{a})\right)=\left(E_{k}|_{k}\sigma_{\mathfrak{b}}\right)(\Gamma,\overline{\nu},z,\mathfrak{a}).
\]

Finally, since $\xi_{2-k}(\mathcal{E}_{2-k})$ is holomorphic, by using \eqref{eq 2.2} we have
\[
\Delta_{2-k}(\mathcal{E}_{2-k})=-\xi_k\left(\xi_{2-k}(\mathcal{E}_{2-k})\right)=0.
\] 
\end{proof}
\end{thm}

From Theorem \ref{thm 4.2}, it is clear that $\mathcal{E}_{2-k}(\Gamma,\nu, z, \mathfrak{a}) \in H_{2-k}^{\#}(\Gamma, \nu)$. Let us denote by $\mathfrak{E}_{2-k}^{\#}(\Gamma,\nu)$, the subspace of $H_{2-k}^{\#}(\Gamma,\nu)$ 
generated by the mock Eisenstein series corresponding to all the cusps of $\Gamma$. That is,
\[
\mathfrak{E}_{2-k}^{\#}(\Gamma,\nu):=\text{Span}_{\mathbb{C}}\{\mathcal{E}_{2-k}(\Gamma,\nu,z,\mathfrak{a})~|~\mathfrak{a}~\text{is a cusp of}~\Gamma\}.
\]
Also, let $\mathfrak{E}_{k}(\Gamma,\overline{\nu})$ denote the Eisenstein space in 
$M_k(\Gamma,\overline{\nu}).$
\begin{cor}\label{cor 4.3}
For $\frac{5}{2}\leq k\in\frac{1}{2}\mathbb{Z},$ the restriction of the shadow operator $\xi_{2-k}$ to $\mathfrak{E}^{\#}_{2-k}(\Gamma,\nu)$ $($which also we denote by $\xi_{2-k})$ is an isomorphism of the vector spaces 
$\mathfrak{E}^{\#}_{2-k}(\Gamma,\nu)$ and $\mathfrak{E}_{k}(\Gamma,\overline{\nu}).$
\begin{proof}
The proof is immediate by using the fact that there are no holomorphic modular forms of negative weight. Indeed, Theorem \ref{thm 4.2} shows that the above restriction of the shadow map is surjective. Since the holomorphic part of any $f\in\mathfrak{E}^{\#}_{2-k}(\Gamma,\nu)$ does not contain any negative power of $q$, the kernel of the restriction of $\xi_{2-k}$ to $\mathfrak{E}^{\#}_{2-k}(\Gamma,\nu)$ is $M_{2-k}(\Gamma,\overline{\nu})$ and hence $\{0\}$ as $2-k\leq -\frac{1}{2}$. 
\end{proof}
\end{cor} 
\begin{cor}\label{cor 4.4}
The space $\mathfrak{E}^{\#}_{2-k}(\Gamma,\nu)$ is finite dimensional. Moreover, for $2<k\in\mathbb{Z}$ we have that 
\[
\mathrm{dim}\mathfrak{E}^{\#}_{2-k}(\Gamma_0(N),\chi)=\sum_{\substack{C \mid N\\\text{gcd}(C, N / C) \mid N / m_{\overline{\chi}}}} \phi(\operatorname{gcd}(C, N / C)),
\] 
where $m_{\overline{\chi}}$ is the conductor of $\overline{\chi}.$
\begin{proof}
This follows from Corollary \ref{cor 4.3} and \cite[Proposition 8.5.15]{CS}.
\end{proof}
\end{cor}

We now define mock Eisenstein series of weight $2-k$ for $k=2,3/2$ and $\Gamma=\Gamma_0(N)$, 
$\nu=\Psi_{2-k,\chi}$. To make the definition, we first prove the following proposition.
\begin{prop}\label{prop:mock01/2anacont}
The nonhomolomrphic Eisenstein series 
$\mathcal{E}_{2-k}(\Gamma_0(N),\Psi_{2-k,\chi},z,\mathfrak{a},s)$ given by \eqref{eq:nonholmock} corresponding to the cusps $\mathfrak{a}=\infty, 0$ can be analaytically 
continued to $s=k-1$ for $k=2$ if $\chi\neq \mathbf{1}_N$, and for $k=3/2$ if $\chi^2\neq \mathbf{1}_N$.
\begin{proof}
The nonholomorphic Eisenstein series $\mathcal{E}_{2-k}(\Gamma_0(N),\Psi_{2-k,\chi},z,\mathfrak{a},s)$ is defined by \eqref{eq:nonholmock} for Re$(s)>k/2$. 
By using \eqref{eq 4.7} and \eqref{eq:hatetas}, we see that 
\begin{equation}\label{eq:mock2}
\mathcal{E}_{2-k}(\Gamma_0(N),\Psi_{2-k,\chi},z,\mathfrak{a},s)=\frac{\delta_{\mathfrak{a}\infty}}{s}y^s+\sum_{n\in\mathbb{Z}}\widehat{\eta}_{\mathfrak{a}\infty}(n,s,y)q^n,
\end{equation}
where 
\[
\widehat{\eta}_{\mathfrak{a}\infty}(n,s,y)=\frac{y^s\varrho^{2-k}_n(s,y)}{s}\sum_{c>0}\frac{S_{\mathfrak{a}\infty}(\Psi_{2-k,\chi},n,c)}{c^{2(1+s)-k}}.
\]
From \eqref{eq:rho-}, \eqref{eq:rho+} and \eqref{eq:rho0}, we see that the factor 
$\frac{y^s\varrho^{2-k}_n(s,y)}{s}$ is well defined at $s=k-1$ for $k=2$ and $3/2$. 

From the calculations done in Section \ref{subsec:intweighteis}, we see that the sum corresponding to $k=2$,
\[
\sum_{c>0}\frac{S_{\mathfrak{a}\infty}(\Psi_{0,\chi},n,c)}{c^{2s}}
\]
is finite for all $n\neq 0$ at $s=1$. In the paragraph before Definition \ref{def:weight2}, 
we noted that the above series is 
finite for $n=0$ at $s=1$ if $\chi \neq \mathbf{1}_N$. 

For $\frac{3}{2} \leq k\in\mathbb{Z}+\frac{1}{2}$ and Re$(s)>k/2$, from Section \ref{subsec 4.2} we have
\[
\sum_{c>0}\frac{S_{\infty\infty}(\Psi_{2-k,\chi},n,c)}{c^{2(1+s)-k}}=A_{2-k}(n,\chi,s),\quad \sum_{c>0}\frac{S_{0\infty}(\Psi_{2-k,\chi},n,c)}{c^{2(1+s)-k}}=B_{2-k}(n,\chi,s).
\]  
From Theorem \ref{thm:anacontAB}, we have meromorphic continuations of $A_{2-k}(n,\chi,s)$ and 
$B_{2-k}(n,\chi,s)$ to the whole $s$-plane, for all $n$. From the discussion before Definition \ref{def:eis3/2anacont}, we see that $A_{1/2}(n,\chi,s)$ and $B_{1/2}(n,\chi,s)$ 
are finite at $s=1/2$ if $\chi^2\neq \mathbf{1}_N$.   
\end{proof} 
\end{prop}

\begin{Def}\label{def:mockeis01/2}
\normalfont
We define mock Eisenstein series $\mathcal{E}_{2-k}(\Gamma_0(N),\Psi_{2-k,\chi},z,\mathfrak{a},s)$ of weight 
$2-k$, where $k=2, 3/2$ corresponding to the cusps $\mathfrak{a}=\infty,0$ by evaluating the corresponding analytic continuation of \eqref{eq:mock2} at $s=k-1$. Note that this definition is made under the assumption that $\chi\neq \mathbf{1}_N$ if $k=2$, and $\chi^2\neq \mathbf{1}_N$ if $k=3/2$. 
\end{Def}
\begin{remark}
Proposition \ref{prop:mock01/2anacont} and \eqref{eq 4.10} ensure that the above defined mock Eisenstein series satisfy correct modular transformations. Moreover, we see that these mock Eisenstein series have similar Fourier series expansions as in \eqref{eq 4.11} with analytically continued Fourier coefficients described in Proposition \ref{prop:mock01/2anacont}. Thus 
we conclude that these mock Eisenstein series are harmonic Maass forms of appropriate growth at the cusps 
(as in Theorem \ref{thm 4.2}) with shadows being the holomorphic Eisenstein series of weight $k=2, 3/2$ and character $\overline{\chi}$.\\
{\textbf{Notation:}} For $\frac{3}{2}\leq k\in \frac{1}{2}\mathbb{Z}$, 
we denote the mock Eisenstein series $\mathcal{E}_{2-k}(\Gamma_0(N),\Psi_{2-k,\chi},z,\infty)$ and 
$\mathcal{E}_{2-k}(\Gamma_0(N),\Psi_{2-k,\chi},z,0)$ corresponding to the cusps $\infty$ and $0$ by 
$\mathcal{E}_{2-k}(\Gamma_0(N),\chi,z)$ and $\mathcal{F}_{2-k}(\Gamma_0(N),\chi,z)$, respectively.
\end{remark}

\section{Proof of Theorems \ref{thm 1.1} and \ref{thm 1.2}}
\label{sec 6}
\noindent As proved in Theorem \ref{thm 4.2}, the shadows of the mock Eisenstein series $\mathcal{E}_{2-k}(\Gamma_0(N),\chi,z)$ and $\mathcal{F}_{2-k}(\Gamma_0(N),\chi,z)$ are the Eisenstein series  $E_{k}(\Gamma_0(N),\overline{\chi},z)$ and $F_{k}(\Gamma_0(N),\overline{\chi},z)$ respectively. From \eqref{eq 2.3} we see that the coefficients of the nonholomorphic parts of these mock Eisenstein series are directly related to the coefficients of its shadows, which have been
already calculated in Sections \ref{sec 3} and \ref{sec 4} for $k \geq 5/2$ and $k =2, 3/2$, respectively. Now we compute the coefficients of the holomorphic parts $\mathcal{E}_{2-k}^+(\Gamma_0(N),\chi,z)$ and $\mathcal{F}_{2-k}^+(\Gamma_0(N),\chi,z)$ of these mock Eisenstein series.

\subsection{Proof of Theorem \ref{thm 1.1}}
For any $2 \leq k \in \mathbb{Z}$, we write 
\[
\begin{split}
\mathcal{E}_{2-k}^+(\Gamma_0(N),\chi,z)=\sum_{n=0}^{\infty}a_{2-k}^+(n,\chi)q^n,\\ \mathcal{F}_{2-k}^+(\Gamma_0(N),\chi,z)=\sum_{n=0}^{\infty}b_{2-k}^+(n,\chi)q^n.
\end{split}
\]
From \eqref{eq 4.12} and proof of Proposition \ref{prop:mock01/2anacont}, we have 
\[
a_{2-k}^+(n,\chi)=\frac{\pi}{k-1}(-2i)^{2-k}\sum\limits_{c=1}^{\infty}\frac{S_{\infty\infty}(\Psi_{2-k},n,c)}{c^{k}}.
\]
By following the computations done in Section \ref{subsec:intweighteis}, we have
\[
a_{2-k}^+(n,\chi)=\frac{(-2i)^{2-k}\pi n^{1-k}}{N^{k}(k-1)}\tau(\overline{\chi}^0)\sigma_{k-1}^{\overline{\chi}}(n),~~~n\neq 0.
\] 
Next, by using \eqref{eq:Sinfinfintweight} we have
\[
S_{\infty\infty}(\Psi_{2-k},0,c)=\begin{cases}
\sideset{}{^{*}}\sum\limits_{m\bmod c}\overline{\chi}(m)&\text{if $N|c$}\\0&\text{otherwise}.
\end{cases}
\]
So we have
\[
\begin{split}
a_{2-k}^+(0,\chi)&=\frac{(-2i)^{2-k}\pi}{(k-1)}\sum_{c=1}^{\infty}\left(\sideset{}{^{*}}\sum\limits_{m\bmod Nc}\overline{\chi}(m)\right)(Nc)^{-k}\\&=\frac{(-2i)^{2-k}\pi}{(k-1)N^k}\sum_{c=1}^{\infty}\left(\sideset{}{^{*}}\sum\limits_{m\bmod c}\overline{\chi}(m)\right)c^{1-k},
\end{split}
\]
here we have used \cite[Eq. 7.2.46]{miyake} in the second step.
We now use the fact that 
\[
\sideset{}{^{*}}\sum\limits_{m\bmod c}\overline{\chi}(m)=\begin{cases}\phi(N)&\text{if $\chi=\mathbf{1}_N$}\\0&\text{if $\chi\neq\mathbf{1}_N$}\end{cases}=\begin{cases}N\prod\limits_{p|N}(1-p^{-1})&\text{if $\chi=\mathbf{1}_N$}\\0&\text{if $\chi\neq\mathbf{1}_N$},
\end{cases}
\]
to get
\begin{equation}\label{eq:mockak}
a_{2-k}^+(0,\chi)=\frac{(-2i)^{2-k}\pi}{N^{k-1}(k-1)}\zeta(k-1) \times \begin{cases}\prod\limits_{p|N}(1-p^{-1})&\text{if $\chi=\mathbf{1}_N$}\\0&\text{if $\chi\neq\mathbf{1}_N$}.
\end{cases}
\end{equation}
Similarly, we have
\[
b_{2-k}^+(n,\chi)=\frac{\pi}{k-1}(-2i)^{2-k}\sum\limits_{c=1}^{\infty}\frac{S_{0\infty}(\Psi_{2-k},n,c)}{c^{k}}.
\]
By following the computations done in Section \ref{subsec:intweighteis}, we get 
\[
b_{2-k}^+(n,\chi)=\frac{(-2i)^{2-k}\pi n^{1-k}L(k,\chi)^{-1}}{N^{k/2}(k-1)}\widetilde{\sigma}_{k-1}^{\chi}(n),~~~n\neq 0.
\]
Finally using \eqref{eq:S0infintweit}, we have 
\[
\begin{split}
b_{2-k}^+(0,\chi)&=\frac{(-2i)^{2-k}\pi}{k-1}\sum_{\substack{c=1\\\text{gcd}(c,N)=1}}^{\infty}(c\sqrt{N})^{-k}\chi(-c)\phi(c)\\&=\frac{(2i)^{2-k}\pi}{N^{k/2}(k-1)}\sum_{\substack{c=1\\\text{gcd}(c,N)=1}}^{\infty}c^{-k}\chi(c)\phi(c),
\end{split}
\]
here we used $\chi(-1)=(-1)^k$.
We can easily see that
\[
\sum_{\substack{c=1\\\text{gcd}(c,N)=1}}^{\infty}c^{-k}\chi(c)\phi(c)=\frac{L(k-1,\chi)}{L(k,\chi)},
\]
Since $\chi\neq\textbf{1}_N$ if $k=2$, the above expression is finite. 
This gives us
\begin{equation}\label{eq:mockbk}
b_{2-k}^+(0,\chi)=\frac{(2i)^{2-k}\pi L(k,\chi)^{-1}}{N^{k/2}(k-1)}L(k-1,\chi).
\end{equation}
\hfill$\square$

\subsection{Proof of Theorem \ref{thm 1.2}} 
For any $3/2 \leq k \in \mathbb{Z}+\frac{1}{2}$, we write
\begin{equation}
\begin{split}
\mathcal{E}_{2-k}^+(\Gamma_0(N),\chi,z)=\sum_{n=0}^{\infty}A_{2-k}^+(n,\chi)q^n,\\ \mathcal{F}_{2-k}^+(\Gamma_0(N),\chi,z)=\sum_{n=0}^{\infty}B_{2-k}^+(n,\chi)q^n.
\end{split}
\label{eq:E+F+series}
\end{equation}
By following the computations done in Section \ref{subsect:half-integral} and \eqref{eq 4.12}, 
for $k\geq \frac{5}{2}$ we have   
\begin{equation}
\begin{split}
& A_{2-k}^+(n,\chi)=\frac{(-2i)^{2-k}\pi}{(k-1)}\sum_{\substack{c=1\\N|c}}^{\infty}\left(\sum_{m \bmod ~c}\overline{\chi}(m)\left(\frac{c}{m}\right)\varepsilon_{m}^{2(2-k)}e^{2\pi i\frac{nm}{c}}\right)c^{-k},\\&
B_{2-k}^+(n,\chi)=\frac{(-2i)^{2-k}i^{2(2-k)}\pi}{N^{k/2}(k-1)}\sum_{\substack{c=1\\\text{gcd}(c,N)=1}}^{\infty}\left(\chi(c)\varepsilon_{c}^{-2(2-k)}\sum_{m \bmod ~c}\left(\frac{m}{c}\right)e^{2\pi i\frac{nm}{c}}\right)c^{-k}. 
\end{split}
\label{eq:E+F+coeff}
\end{equation}
Comparing with \eqref{eq:coeffnonholeisboth}, we have 
\begin{equation}
A_{2-k}^+(n,\chi)=\frac{(-2i)^{2-k}\pi}{(k-1)}A_{2-k}(n,\chi,k-1),\ B_{2-k}^+(n,\chi)=\frac{(2i)^{2-k}\pi}{N^{k/2}(k-1)}B_{2-k}(n,\chi,k-1).
\label{eq:mockeishalfint}
\end{equation}
From Proposition \ref{prop:mock01/2anacont}, we see that \eqref{eq:mockeishalfint} holds true for $k=3/2$ as well with the understanding that in this case the right hand side of \eqref{eq:mockeishalfint} is the value of the analytic continuation of these coefficients. The claimed expression of $B^+_{2-k}(n,\chi)$ now follows 
from Theorem \ref{thm:anacontAB}. Also, we have 
\[
\begin{split}
A_{2-k}^+(n,\chi)&=\frac{(-2i)^{2-k}\pi}{(k-1)}B_{2-k}(n,\overline{\chi},k-1)C_{2-k}(n,\overline{\chi}, k-1)\\&=(-iN)^{k/2}B_{2-k}^+(n,\chi)C_{2-k}(n,\overline{\chi}, k-1).
\end{split}
\]
Since $\chi^2\neq\textbf{1}_N$ if $k=\frac{3}{2}$, the discussion before Definition \ref{def:eis3/2anacont} implies that the coefficients are well-defined for $\frac{3}{2}\leq k\in\mathbb{Z}+\frac{1}{2}$. \hfill$\square$

\section{Proof of Theorems \ref{thm 1.3} and \ref{thm 1.4}}
\label{sec 7}
\noindent Throughout this section, we denote the trivial character $\mathbf{1}_4$ modulo $4$ also by $1$ 
and the usage will be clear from the context.
To prove Theorems \ref{thm 1.3} and \ref{thm 1.4}, we need the following proposition. 
\begin{prop}\label{prop:thetapower}
We have that 
\[
\begin{split}
&\Theta^3(z)=E_{3/2}(\Gamma_0(4),1,z)+\left(\frac{1-i}{\sqrt{2}}\right)^3F_{3/2}(\Gamma_0(4),1,z)\\
&\Theta^5(z)=E_{5/2}(\Gamma_0(4),1,z)+\left(\frac{1-i}{\sqrt{2}}\right)^5F_{5/2}(\Gamma_0(4),1,z)\\
&\Theta^7(z)=E_{7/2}(\Gamma_0(4),1,z)+\left(\frac{1-i}{\sqrt{2}}\right)^7F_{7/2}(\Gamma_0(4),1,z).
\end{split}
\]
\begin{proof}
We begin by calculating the value of $\Theta$ at the cusps $\infty$ and $0$. Obviously, at the cusp $\infty$ the value is 1. For the cusp $0$, the scaling matrix for $\Gamma_0(4)$ is given by 
\[
\sigma_{0}=\begin{pmatrix}
0&-1/\sqrt{4}\\\sqrt{4}&0
\end{pmatrix}.
\]
We have that 
\[
(\Theta |_{1/2}\sigma_0)(z)=(\sqrt{4}z)^{-1/2}\Theta(-1/4z)=\sqrt{2}(4z)^{-1/2}\Theta(-1/4z)=\sqrt{2}(\Theta|_{1/2}S)(4z),~~~~S=\begin{pmatrix}
0&-1\\1&0
\end{pmatrix}.
\]
We know that (cf. Exercise 7, \S IV.1 of \cite{NK}) 
\[
\lim_{z\to \infty}(\Theta|_{1/2}S)(4z)=\frac{1-i}{2}.
\]
So the value of $\Theta$ at the cusp $0$ is $(1-i)/\sqrt{2}.$ Now noting that $\Theta^{\frac{\kappa}{2}}|_{\kappa/2}S=(\Theta|_{1/2}S)^{\kappa}$, we see that 
the last two equalities of the proposition follow from the orthogonal decomposition of the space of modular forms into the Eisenstein space and cusp forms \cite{He,Pe}, and the fact that each of the spaces of cusp forms of weights $5/2$ and $7/2$ are trivial. 

Now we prove the first relation by using the Fourier expansions of $E_{3/2}$ and $F_{3/2}$ in \eqref{eq:eisnhEF}. 
By using \eqref{eq:coeffnonholeisboth}, we have
\[
\begin{split}
A_{3/2}(n,1,s)&=\sum_{\substack{c=1\\4|c}}^{\infty}\left(\sum_{m\bmod~c}\varepsilon_m^{3}\left(\frac{c}{m}\right)e^{2\pi i\frac{mn}{c}}\right)c^{-3/2-2s}\\&=2^{-3/2-2s}\sum_{\substack{c=2\\c~\text{even}}}^{\infty}\left(\sum_{m\bmod~2c}\varepsilon_m^{3}\left(\frac{2c}{m}\right)e^{\pi i\frac{mn}{c}}\right)c^{-3/2-2s}\\&=2^{-3/2-2s}\sum_{\substack{c=2\\c~\text{even}}}^{\infty}\left(\sum_{m\bmod~2c}e^{\frac{\pi i(m-1)}{4}}\left(\frac{c}{m}\right)e^{\pi i\frac{mn}{c}}\right)c^{-3/2-2s},
\end{split}
\]
here we have used the identity
\begin{equation}
\varepsilon_m^{3}\left(\frac{2}{m}\right)=e^{\frac{\pi i(m-1)}{4}}.
\label{eq:iden2/m}
\end{equation}
Similarly, we have
\[
\begin{split}
B_{3/2}(n,1,s)&=\sum_{\substack{c=1\\\text{gcd}(4,c)=1}}^{\infty}\varepsilon_c^{-3}\left(\sum_{m\bmod~c}\left(\frac{m}{c}\right)e^{2\pi i\frac{mn}{c}}\right)c^{-3/2-2s}\\&=\sum_{\substack{c=1\\c~\text{odd}}}^{\infty}\varepsilon_c^{-3}\left(\frac{2}{c}\right)\left(\sum_{m\bmod~c}\left(\frac{2m}{c}\right)e^{2\pi i\frac{mn}{c}}\right)c^{-3/2-2s}\\&=\sum_{\substack{c=1\\c~\text{odd}}}^{\infty}e^{\frac{\pi i(1-c)}{4}}\left(\sum_{m\bmod~2c}\left(\frac{m}{c}\right)e^{\pi i\frac{mn}{c}}\right)c^{-3/2-2s},
\end{split}
\]
where we have again used \eqref{eq:iden2/m} and changed the summation variable of the inner sum in the second step. Now, we have 
\[
\begin{split}
E_{3/2}(\Gamma_0(4),1,z,s)+&\left(\frac{1-i}{\sqrt{2}}\right)^3F_{3/2}(\Gamma_0(4),1,z,s)\\&=1+\sum_{n\in\mathbb{Z}}\left(A_{3/2}(n,1,s)+\frac{i^3}{2^{3/2+2s}}\left(\frac{1-i}{\sqrt{2}}\right)^3B_{3/2}(n,1,s)\right)\varrho^k_n(s,y)q^n.
\end{split}
\]
Using above computations, we have
\[
\begin{split}
A_{3/2}(n,1,s)+\frac{i^{3}}{2^{3/2+2s}}&\left(\frac{1-i}{\sqrt{2}}\right)^3B_{3/2}(n,1,s)\\&=\frac{1}{2^{3/2+2s}}\left[\sum_{\substack{c=2\\c~\text{even}}}^{\infty}\left(\sum_{m\bmod~2c}e^{\frac{3\pi i(1-m)}{4}}\left(\frac{c}{m}\right)e^{\pi i\frac{mn}{c}}\right)c^{-3/2-2s}\right.\\&\left.\hspace{3cm}+\sum_{\substack{c=1\\c~\text{odd}}}^{\infty}e^{\frac{3\pi ic}{4}}\left(\sum_{m\bmod~2c}\left(\frac{m}{c}\right)e^{\pi i\frac{mn}{c}}\right)c^{-3/2-2s}\right],
\end{split}
\]
where we have used the identities 
\begin{equation}
e^{\frac{\pi i(m-1)}{4}}=e^{\frac{3\pi i(1-m)}{4}},\quad e^{\frac{\pi i(1-c)}{4}}=e^{\frac{3\pi i(c-1)}{4}},\quad m,c~\text{odd};\quad \left(\frac{1-i}{\sqrt{2}}\right)^3=e^{-\frac{3\pi i}{4}}. 
\label{eq:idenexp}
\end{equation}

Introduce the Dirichlet series  
\begin{equation}
Z_n(s):=\sum_{c=1}^{\infty}\frac{\mathcal{S}(n,c)}{c^{\frac{1}{2}+s}},
\label{eq 5.1}
\end{equation}
where
\[
\mathcal{S}(n,c):=\sum_{m\bmod 2c}\lambda(m, c) e^{\pi i\frac{mn}{c}}
\]
with 
\begin{equation*}
\lambda(m, c):=\begin{cases}
e^{\frac{3\pi ic}{4}}\left(\frac{m}{c}\right) & \text {if $c$ is odd, $m$ is even} \\
e^{\frac{3\pi i(1-m)}{4}}\left(\frac{c}{m}\right) & \text {if  $c$ is even, $m$ is odd} \\
0 & \text {otherwise}.
\end{cases}
\label{eq 5.2}
\end{equation*}
With this notation, we have
\begin{equation}
E_{3/2}(\Gamma_0(4),1,z,s)+\left(\frac{1-i}{\sqrt{2}}\right)^3F_{3/2}(\Gamma_0(4),1,z,s)=1+\frac{1}{2^{3/2+2s}}\sum_{n\in\mathbb{Z}}Z_n(2s+1)\varrho^k_n(s,y)q^n.
\label{eq:EFtheta}
\end{equation}

The Dirichlet series $Z_n(2s+1)$ admits an analytic continuation (see \cite[Section 3.1, Section 3.2]{RW}) which is finite at $s=0$ for every $n\in\mathbb{Z}$. Using this and \eqref{eq 4.4}, we see that the right hand side of \eqref{eq:EFtheta} at $s=0$ is equal to 
\[
1+2\pi i^{-3/2}\sum_{n=1}^{\infty}Z_n(1)n^{1/2}q^n=\Theta^3.
\]
Here we have used the fact that the Fourier coefficient $r_3(n)$ of $\Theta^3$ is given by \cite[Eq. (2.6)]{RW}  
\[
r_3(n)=2i^{-3/2}\pi n^{1/2}Z_n(1).
\]

Finally, we need to show that the value of the analytic continuation of the left hand side of \eqref{eq:EFtheta} 
at $s=0$ is equal to
\[
E_{3/2}(\Gamma_0(4),1,z)+\left(\frac{1-i}{\sqrt{2}}\right)^3F_{3/2}(\Gamma_0(4),1,z).
\]  
Indeed we have
\[
\begin{split}
E_{3/2}(&\Gamma_0(4),1,z,0)+\left(\frac{1-i}{\sqrt{2}}\right)^3F_{3/2}(\Gamma_0(4),1,z,0)\\&=1+\sum_{n\in\mathbb{Z}}\left(A_{3/2}(n,1,0)\varrho^k_n(0,y)+\frac{i^3}{2^{3/2}}\left(\frac{1-i}{\sqrt{2}}\right)^3B_{3/2}(n,1,0)\varrho^k_n(0,y)\right)q^n\\&=1+\sum_{n\in\mathbb{Z}}A_{3/2}(n,1,0)\varrho^k_n(0,y)q^n+\frac{i^3}{2^{3/2}}\left(\frac{1-i}{\sqrt{2}}\right)^3\sum_{n\in\mathbb{Z}}B_{3/2}(n,1,0)\varrho^k_n(0,y)q^n\\&=E_{3/2}(\Gamma_0(4),1,z)+\left(\frac{1-i}{\sqrt{2}}\right)^3F_{3/2}(\Gamma_0(4),1,z).
\end{split}
\]
We need to justify the second last step which follows if we can show that the series defining $E_{3/2}(\Gamma_0(4),1,z)$ and $F_{3/2}(\Gamma_0(4),1,z)$ converge. 
Now we show the absolute convergence of the series. We separate the series into two parts $n> 0$ and $n<0$.
\[
E_{3/2}(\Gamma_0(4),1,z)=1+A_{3/2}(0,1,0)\varrho^k_0(0,y)+\sum_{n>0}A_{3/2}(n,1,0)\varrho^k_n(0,y)q^n+\sum_{n<0}A_{3/2}(n,1,0)\varrho^k_n(0,y)q^n.
\]
By Theorem \ref{thm 4.1}, we see that 
\[
A_{3/2}(n,1,0)=B_{3/2}(n,1,0)C_{3/2}(n,1,0),
\]
where $C_{3/2}$ is a finite Dirichlet series for $n\neq 0$. Moreover
\[
|C_{3/2}(0,1,0)|\leq \sum_{4|M|4^{\infty}}M^{-1/2}< \infty.
\]
The coefficient $B_{3/2}(n,1,0)$ has a pole due to $L(1,\omega_n)$ when $n\leq 0$ and $-n$ is a square. This pole is cancelled by $\Gamma(0)^{-1}$ factor in $\varrho^k_n(0,y)$. Using \eqref{eq 4.4} and \eqref{eq 4.6}, we see that 
\[
|A_{3/2}(0,1,0)\varrho^k_0(0,y)| < \infty.
\]
For $n>0$, we also have 
\[
A_{3/2}(n,1,0)\varrho^k_n(0,y)=O(n^{\alpha})
\]
for some $\alpha>0$. Finally for $n<0$, we see that 
\[
A_{3/2}(n,1,0)\varrho^k_n(0,y)=a(n)y^{-3/2}e^{4\pi ny}\Omega(-4\pi ny,-1/2,1),
\]
where $a(n)=O(n^{\nu})$ for some $\nu>0$. Finally using \eqref{eq:omeinga}, we have 
\[
A_{3/2}(n,1,0)\varrho^k_n(0,y)=a'(n)\Gamma(1/2,-4\pi ny),
\]
where $a'(n)=O(n^{\nu'})$ for some $\nu'>0$. Thus we obtain 
\[
E_{3/2}(\Gamma_0(4),1,z)=1+a^-(0)y^{-1/2}+\sum_{n>0}a^{+}(n)q^n+\sum_{n<0}a^{-}(n)\Gamma(1/2,-4\pi ny)q^n,
\]
where $a^{\pm}(n)=O(|n|^{\sigma})$ for some $\sigma>0.$
The convergence is now immediate  by using the asymptotic relation \cite[Eq. (4.6)]{Ono1}
\[
\Gamma(s,x)\sim x^{s-1}e^{-x},\quad x\in\mathbb{R},\quad |x|\to \infty.
\]
The convergence of $F_{3/2}(\Gamma_0(4),1,z)$ can be proved similarly.
\end{proof}
\end{prop} 
The following result is known but we could not find it written down explicitly in existing literature. 
We therefore write it as a lemma and provide a sketch of the proof below.
\begin{lemma}\label{lemma:thetaeigenform}
Let $p$ be an odd prime and $\kappa\in\{3,5,7\}$. Then  $\Theta^{\kappa}$ are Hecke eigenforms of all the Hecke operators $T(p^2)$ with eigenvalue $1+p^{\kappa}.$
\begin{proof}
By using the Shimura lifting given in \cite[Section 1, Page 78]{ash}, we see that the images of $\Theta^5(z)$ and $\Theta^7(z)$ under the first Shimura map are $\frac{1}{24}E_4(z)$ and $-\frac{1}{6}L(-2, \chi_{-1}) \left(E_6(z)+ 8E_6(2z)\right)$ respectively, where $\chi_{-1}= \left(\frac{-1}{\cdot}\right)$. Therefore, $\Theta^5$ and $\Theta^7$ are Hecke eigenforms  with respect to $T_{p^2}$, $p$ odd, with eigenvalues $1+ p^3$ and $1+ p^5$ respectively. 

Also, $M_{3/2}(\Gamma_0(4))$ is one dimensional, therefore $T_{p^2}\Theta^3= \lambda_p\Theta^3$ for any odd prime $p$. Comparing the constant coefficient on both sides, we get $\lambda_p= 1+ p$.
\end{proof} 
\end{lemma} 

\subsection{Proof of Theorem \ref{thm 1.3}} 

For $\kappa\in\{5,7\}$, consider the functions
\[
F_{\Theta^{\kappa}}(z)=\mathcal{E}_{2-\kappa/2}(\Gamma_0(4),1,z)+\left(\frac{1+i}{\sqrt{2}}\right)^{\kappa}\mathcal{F}_{2-\kappa/2}(\Gamma_0(4),1,z).
\]
By using Theorem \ref{thm 4.2} and Proposition \ref{prop:thetapower}, we get that $F_{\Theta^{\kappa}}$ is a harmonic Maass form with shadow $\Theta^{\kappa}$ for any $\kappa\in\{5,7\}$. 
Now we compute the Fourier coefficients $c_{F_{\Theta^{\kappa}}}^+(n)$ of the holomorphic part $F_{\Theta^{\kappa}}^+$ 

By using \eqref{eq:E+F+series} and \eqref{eq:E+F+coeff}, we have  
\[
\begin{split}
c_{F_{\Theta^{\kappa}}}^+(n)&=\frac{(-2i)^{2-\kappa/2}\pi}{(\kappa/2-1)}\sum_{\substack{c=1\\4|c}}^{\infty}\left(\sum_{m (\bmod ~c)}\left(\frac{c}{m}\right)\varepsilon_{m}^{2(2-\kappa/2)}e^{2\pi i\frac{nm}{c}}\right)c^{-\kappa/2}\\&+\left(\frac{1+i}{\sqrt{2}}\right)^{\kappa}\frac{(-2i)^{2-\kappa/2}i^{2(2-\kappa/2)}\pi}{4^{\kappa/4}(\kappa/2-1)}\sum_{\substack{c=1\\\text{gcd}(c,4)=1}}^{\infty}\left(\varepsilon_{c}^{-2(2-\kappa/2)}\sum_{m (\bmod ~c)}\left(\frac{m}{c}\right)e^{2\pi i\frac{nm}{c}}\right)c^{-\kappa/2}\\&=\frac{(-2i)^{2-\kappa/2}\pi}{(\kappa/2-1)2^{\kappa/2}}\sum_{\substack{c=1\\c~\text{even}}}^{\infty}\left(\sum_{m (\bmod ~2c)}\left(\frac{c}{m}\right)\left(\frac{2}{m}\right)\varepsilon_{m}^{-\kappa}e^{\pi i\frac{nm}{c}}\right)c^{-\kappa/2}\\&+e^{\frac{\pi i\kappa}{4}}\frac{(-2i)^{2-\kappa/2}i^{-\kappa}\pi}{2^{\kappa/2}(\kappa/2-1)}\sum_{\substack{c=1\\c~\text{odd}}}^{\infty}\left(\varepsilon_{c}^{\kappa}\left(\frac{2}{c}\right)\sum_{m (\bmod ~2c)}\left(\frac{m}{c}\right)e^{\pi i\frac{nm}{c}}\right)c^{-\kappa/2},
\end{split}
\]
where we have used similar manipulations as in the second step of the proof of Proposition \ref{prop:thetapower}. Now we use the following identities.
\begin{equation}
\varepsilon_m^{-2\ell-1}\left(\frac{2}{m}\right)=e^{\frac{\pi i(1-m)}{4}},\quad \varepsilon_c^{2\ell'+1}\left(\frac{2}{c}\right)=e^{\frac{\pi i(1-c)}{4}},\quad \text{$\ell$ odd, $\ell'$ even.}
\label{eq:iden2/mex}
\end{equation}
Using \eqref{eq:idenexp} and \eqref{eq:iden2/mex}, we get
\[
\begin{split}
c_{F_{\Theta^{\kappa}}}^+(n)&=\frac{(-2i)^{2-\kappa/2}\pi}{(\kappa/2-1)2^{\kappa/2}}\sum_{\substack{c=1\\c~\text{even}}}^{\infty}\left(\sum_{m (\bmod ~2c)}e^{(-1)^{\frac{\kappa+1}{2}}\frac{3\pi i(m-1)}{4}}\left(\frac{c}{m}\right)e^{\pi i\frac{nm}{c}}\right)c^{-\kappa/2}\\&+e^{-\frac{\pi i\kappa}{4}}\frac{(-2i)^{2-\kappa/2}\pi}{2^{\kappa/2}(\kappa/2-1)}\sum_{\substack{c=1\\c~\text{odd}}}^{\infty}\left(\sum_{m (\bmod ~2c)}e^{(-1)^{\frac{\kappa+1}{2}}\frac{3\pi i(1-c)}{4}}\left(\frac{m}{c}\right)e^{\pi i\frac{nm}{c}}\right)c^{-\kappa/2}.
\end{split}
\]

We thus have
\begin{equation*}
c_{F_{\Theta^{\kappa}}}^+(n)=\frac{(-2i)^{2-\kappa/2}\pi}{(\kappa/2-1)2^{\kappa/2}}.\begin{cases}
Z_n\left(\frac{\kappa-1}{2}\right)&\kappa=5\\ & \\\overline{Z_{-n}\left(\frac{\kappa-1}{2}\right)}-2\overline{Z_{-n}^{\text{odd}}(\frac{\kappa}{2})}&\kappa=7,
\end{cases}
\end{equation*} 
where $Z_n(s)$ is as in \eqref{eq 5.1} and $Z_{n}^{\text{odd}}(s)$ is given by 
\[
Z_{n}^{\text{odd}}(s)=\sum_{\substack{c=1\\c~\text{odd}}}^{\infty}\frac{\mathcal{S}(n,c)}{c^{s}}.
\]

The analytic continuations of $Z_n(s)$ and $Z_{n}^{\text{odd}}(s)$ respectively have been worked out in \cite{RW}. By using \cite[Theorems 3.1, 3.2]{RW}, we have the following. Write $n=f^{2} d \neq 0$ with $d$ square-free and $f=2^{q} w$ with $w$ odd, then
\begin{equation*}
Z_{n}^{\text{odd}}(s)=e^{\frac{3 \pi i}{4}} \frac{L\left(s, \psi_{-n}\right)}{\zeta(2 s)} w^{1-2 s} T_{s}^{\psi_{-n}}(w) \frac{1-\psi_{-n}(2) 2^{-s}}{1-2^{-2 s}}
\label{eq:zodd}
\end{equation*}
and 
\begin{equation}
Z_{n}(s)=Z_{n}^{\text {odd}}(s) R_{n}(s),
\label{eq:zns}
\end{equation}
where 
$$
R_{n}(s):=1+2^{-s}-2^{1-s} R_{n}^{*}(s)
$$
with
\begin{equation*}
R_{n}^{*}(s):=\begin{cases}
0 & \text { if } n \equiv 1,2(\bmod~ 4) \\
\frac{1-2^{-2 s}}{1-\psi_{-n}(2) 2^{-s}} 2^{Q(1-2 s)} T_{s}^{\psi_{-n}}\left(2^{Q}\right) & \text { otherwise, }
\end{cases}
\label{eq:rnstar}
\end{equation*}
and $Q$ is defined in \eqref{eq:Q}. Moreover 
\begin{equation*}
Z_{0}(s)=e^{\frac{3 \pi i}{4}} \frac{\zeta(2 s-1)}{\zeta(2 s)} \frac{1-2^{-(2 s-1)}-2^{-s}}{1-2^{-2 s}}.
\end{equation*}
The coefficients $c_{F_{\Theta^{\kappa}}}^+(n)$ for $\kappa=5,7$ depend on $Z_n(s)$ with $s=2,3$. 
To evaluate $Z_n(s)$, we 
consider the cases when $-n$ is a square or not. When $-n$ is not a square then  $\psi_{-n}$ is a nontrivial character. This gives
\[
Z_n(s)=e^{\frac{3 \pi i}{4}} \frac{L(s,\psi_{-2})}{\zeta(2 s)} w^{1-2 s} T_{s}^{\psi_{1}}(w) c_n(s),
\]
where 
\[
c_n(s):=\frac{1-\psi_{-n}(2) 2^{-s}}{1+2^{-2s}}R_n(s).
\]
To simplify $c_n(s)$, we compute $T_{s}^{\psi_{-n}}\left(2^{Q}\right)$. Using \eqref{eq 1.3} we have
\[
T_{s}^{\psi_{-n}}\left(2^{Q}\right)=\sum_{\ell=0}^{Q}\mu(2^{\ell})\psi_{-n}(2^{\ell})2^{\ell(s-1)}\sigma_{2s-1}(2^{Q-\ell})=\sigma_{2s-1}(2^Q)-2^{s-1}\psi_{-n}(2)\sigma_{2s-1}(2^{Q-1}),
\]
here we have used the fact that $\mu(2^{\ell})=0$ for $\ell>1$ and $\mu(1)=1,\mu(2)=-1$.
Thus we have 
\begin{equation}
c_n(s)=\begin{cases}\frac{2^s-\psi_{-n}(2)}{2^s-1};\hspace{1cm}n\equiv 1,2(\bmod ~4)&\\&\\\frac{1-\psi_{-n}(2) 2^{-s}}{1+2^{-s}}-2^{1-s+Q(1-2s)}\left(\sigma_{2s-1}(2^Q)-2^{s-1}\psi_{-n}(2)\sigma_{2s-1}(2^{Q-1})\right);&\text{otherwise.}
\end{cases}
\label{eq:cns}
\end{equation}
When $-n$ is a square then $L\left(s, \psi_{-n}\right)=\zeta(s)$, which in addition with \eqref{eq:zns} gives 
\[
Z_n(s)=e^{\frac{3 \pi i}{4}} \frac{\zeta(s)}{\zeta(2 s)} w^{1-2 s} T_{s}^{\psi_{1}}(w) \frac{R_n(s)}{1+2^{-s}}.
\]
Moreover, for $-n$ are square, we have (cf. \cite[Page 1009]{RW})
\[
\frac{R_{n}(s)}{\left(1+2^{-s}\right)}=1-2^{(1-s)}.
\] 
Hence when $-n$ is a square, we have
\begin{equation*}
Z_n(s)=e^{\frac{3 \pi i}{4}} \frac{\zeta(s)}{\zeta(2 s)} w^{1-2 s} T_{s}^{\psi_{1}}(w)(1-2^{1-s}).
\end{equation*}
Thus for $\kappa=5,7$, we get 
\[
Z_n\left(\frac{\kappa-1}{2}\right)=\left\{\begin{array}{ll}
e^{\frac{3 \pi i}{4}} \frac{\zeta(\kappa-2)}{\zeta(\kappa-1)} \frac{1-2^{-(\kappa-2)}-2^{-(\kappa-1)/2}}{1-2^{-\kappa+1}} & \text { if } n=0 \\
e^{\frac{3 \pi i}{4}}\frac{1}{\zeta(\kappa-1)}\frac{T_{(\kappa-1)/2}(w)}{w^{\kappa-2}}\zeta\left(\frac{\kappa-1}{2}\right)\left(1-2^{\frac{3-\kappa}{2}}\right) & \text { if } -n \text { is a square } \\
e^{\frac{3 \pi i}{4}} \frac{L\left(\frac{\kappa-1}{2}, \psi_{-n}\right)}{\zeta(\kappa-1)}  \frac{T_{(\kappa-1)/2}^{\psi_{-n}}(w)}{w^{\kappa-2}} \cdot c_{n}\left(\frac{\kappa-1}{2}\right) & \text { otherwise. }
\end{array}\right.
\]
For $\kappa=7$, we need to evaluate $Z_{0}^{\text{odd}}(3)$. Indeed, by calculations in \cite[Section 3.2]{RW}, we have 
\[
Z_{n}^{\text{odd}}(s)=e^{\frac{3\pi i}{4}}\sum_{\substack{c=1\\c~\text{odd}}}^{\infty}\frac{\phi(c^2)}{c^{2s}}=\frac{\zeta(2(s-1))}{\zeta(2s-1)}.
\]

Finally, we show that $F_{\Theta^{\kappa}}$ are eigenvectors of Hecke operators $T(p^2)$ with eigenvalues $1+p^{2-\kappa}$. Since $\Theta^\kappa$ are Hecke eigenforms with eigenvalues $1+p^{\kappa-2}$ (see Lemma \ref{lemma:thetaeigenform}), by using \eqref{eq 2.5} we have 
\[
F_{\Theta^{\kappa}}|T(p^2)-\left(1+\frac{1}{p^{\kappa-2}}\right)F_{\Theta^{\kappa}}\in M_{2-\kappa/2}(\Gamma_0(4),\Psi_k).
\]
For $\kappa=5,7$, we see that the space $M_{2-\kappa/2}(\Gamma_0(4),\Psi_k)$ is trivial as $2-\kappa/2<0$. 
\hfill$\square$

\subsection{Proof of Theorem \ref{thm 1.4}} 

To get the lift of $\Theta^3$, we imitate the construction employed in the proof of Theorem \ref{thm 1.3} though the case of $\Theta^3$ is a bit more complicated. This time we need to work with the analytic continuations of the involved mock Eisenstein series of weight $1/2$. To this end, we first define the function 
\[
F_{\Theta^3}(z,s)=\mathcal{E}_{1/2}(\Gamma_0(4),\Psi_{1/2},z,\infty,s)+\left(\frac{1+i}{\sqrt{2}}\right)^3\mathcal{E}_{1/2}(\Gamma_0(4),\Psi_{1/2}, z,0,s),
\]
where $\chi=1$ in the multiplier $\Psi_{1/2}$. 
We cannot directly choose $s=1/2$ to get the harmonic Maass form as we did when we defined mock Eisenstein series because it does not lie in the domain of convergence of the right hand side. We will instead find the analytic continuation of the right hand side and then evaluate it at $s=1/2.$ We will show that the analytic continuation $F_{\Theta^3}(z)$ of the function $F_{\Theta^3}(z,s)$ at $s=1/2$ has the form
\[
F_{\Theta^{3}}(z)=\frac{\pi}{2}e^{-\frac{i\pi}{4}}\overline{Z_0(1)}+e^{-\frac{i\pi}{4}}\pi\sum_{n=1}^{\infty}\overline{Z_{-n}(1)}q^n+2\sqrt{y}+e^{-\frac{i\pi}{4}}\sqrt{\pi}\sum_{n<0}\overline{Z_{-n}(1)}\Gamma(1/2,-4\pi ny)q^n.
\]
Using \eqref{eq 4.7} and \eqref{eq:hatetas}, we can easily see that 
\[
F_{\Theta^3}(z,s)=\frac{y^s}{s}\left[1+\sum_{n\in\mathbb{Z}}\left(A_{1/2}(n,1,s)+4^{-1/4-s}\left(\frac{1+i}{\sqrt{2}}\right)^3B_{1/2}(n,1,s)\right)\varrho^{1/2}_n(s,y)q^n\right].
\]
Using \eqref{eq:coeffnonholeisboth} and proceeding as in the proof of Proposition \ref{prop:thetapower}, we have 
\[
\begin{split}
A_{1/2}(n,1,s)+\frac{i}{2^{1/2+2s}}&\left(\frac{1+i}{\sqrt{2}}\right)^3B_{1/2}(n,1,s)\\&=\frac{1}{2^{1/2+2s}}\left[\sum_{\substack{c=2\\c~\text{even}}}^{\infty}\left(\sum_{m\bmod~2c}e^{\frac{3\pi i(m-1)}{4}}\left(\frac{c}{m}\right)e^{\pi i\frac{mn}{c}}\right)c^{-1/2-2s}\right.\\&\left.\hspace{3cm}+\sum_{\substack{c=1\\c~\text{odd}}}^{\infty}e^{-\frac{3\pi ic}{4}}\left(\sum_{m\bmod~2c}\left(\frac{m}{c}\right)e^{\pi i\frac{mn}{c}}\right)c^{-1/2-2s}\right]\\&=\frac{1}{2^{1/2+2s}}\overline{Z_{-n}(2s)}.
\end{split}
\]
Thus we have 
\[
F_{\Theta^3}(z,s)=\frac{y^s}{s}\left[1+\frac{1}{2^{1/2+2s}}\sum_{n\in\mathbb{Z}}\overline{Z_{-n}(2s)}\varrho^{1/2}_n(s,y)q^n\right].
\]
Since $Z_{n}(2s)$ has analytic continuation to $s=1/2$ and $y^s\varrho^{1/2}_n(s,y)$ is also analytic at $s=1/2$, by using \eqref{eq:rho-}, \eqref{eq:rho+} and \eqref{eq:rho0}, we get 
\[
F_{\Theta^{3}}(z)=\frac{\pi}{2}e^{-\frac{i\pi}{4}}\overline{Z_0(1)}+e^{-\frac{i\pi}{4}}\pi\sum_{n=1}^{\infty}\overline{Z_{-n}(1)}q^n+2\sqrt{y}+e^{-\frac{i\pi}{4}}\sqrt{\pi}\sum_{n<0}\overline{Z_{-n}(1)}\Gamma(1/2,-4\pi ny)q^n.
\] 
The explicit value of $Z_{n}(1)$ is calculated in \cite{RW}.
\[
Z_{n}(1)=\begin{cases}
e^{\frac{3 \pi i}{4}} \frac{6}{\pi^{2}} \log (2) & \text { if } n=0 \\
e^{\frac{3 \pi i}{4}} \frac{6}{\pi^{2}} \log (2) \frac{T_{1}^{\psi_1}(w)}{w} & \text { if } -n \text { is a square } \\
e^{\frac{3 \pi i}{4}}\frac{6}{\pi^{2}} L\left(1, \psi_{-n}\right) \frac{T_{1}^{\psi_{-n}}(w)}{w} \cdot c_{n}(1) & \text { otherwise, }
\end{cases}
\]
where $c_n(1)$, computed using \eqref{eq:cns}, is easily seen to be 
\[
c_{n}(1)=\left\{\begin{array}{ll}
2-\psi_{-n}(2) & \text { if } n \equiv 1,2(\bmod~ 4) \\
2^{-Q}\left(1-\psi_{-n}(2)\right) & \text { otherwise. }
\end{array}\right.
\] 
Since $F_{\Theta^3}(z,s)$ satisfies modularity, so does $F_{\Theta^{3}}(z)$ being the analytic continuation of $F_{\Theta^3}(z,s)$. Moreover the growth condition of $F_{\Theta^{3}}(z)$ is clear from the Fourier expansion. Next, the shadow of $F_{\Theta^{3}}(z)$ can be calculated using \eqref{eq 2.3}.
\[
\xi_{1/2}(F_{\Theta^{3}}(z))=1+2\pi i^{-3/2}\sum_{n=1}^{\infty}Z_n(1)n^{1/2}q^n=\Theta^3(z).
\]
The weight $3/2$ hyperbolic Laplacian annihilates $F_{\Theta^{3}}(z)$. Indeed, by using \eqref{eq 2.2} we have
\[
\Delta_{3/2}(F_{\Theta^{3}}(z))=-\xi_{1/2}\left(\xi_{3/2}\left(F_{\Theta^{3}}(z)\right)\right)=-\xi_{1/2}(\Theta^3)=0.
\]

Finally, the fact that $F_{\Theta^{3}}$ is a Hecke eigenform follows by using Lemma \ref{lemma:thetaeigenform} and arguments similar to those provided for $F_{\Theta^{5}}$ and $F_{\Theta^{7}}$.

\section{Proof of Theorem \ref{thm 1.5}}
\label{sec 8}
\noindent To prove Theorem \ref{thm 1.5}, we will employ the construction of Maass-Poincar\'{e} series as in \cite{BO}. To be consistent with the setup of \cite{BO}, we will recall certain aspects of the Fourier expansion at various cusps of $\Gamma_0(N).$ Let $\rho$ be any cusp of $\Gamma_0(N)$. Let $\Gamma_\rho:=\{g \in \Gamma_0(N) | g \rho = \rho\}$ and $\gamma_\rho \in SL_2(\mathbb{Z})$ be such that $\gamma_\rho (\infty)=\rho$. Then 
$\gamma_\rho^{-1} \Gamma_\rho \gamma_\rho$ fixes $\infty$ and hence generated by $-I$ and $\begin{pmatrix} 
1&t_\rho\\0&1\end{pmatrix}$ for some positive integer $t_\rho$. The integer $t_\rho$ is called the {\em width} of the cusp 
$\rho$. Let $g_\rho \in \Gamma_\rho$ such that
$\gamma_\rho^{-1} g_\rho \gamma_\rho = \begin{pmatrix}1&t_\rho\\0&1\end{pmatrix}$. For any $f\in H_k^!(\Gamma_0(N),\chi)$, 
we have
$$
(f|\gamma_\rho)| \begin{pmatrix}1&t_\rho\\0&1\end{pmatrix} = (f| g_\rho)|\gamma_\rho = 
\chi(d_\rho) f|\gamma_\rho, \ g_\rho = \begin{pmatrix} a_\rho & b_\rho\\c_\rho & d_\rho \end{pmatrix}.
$$
Let $\kappa_\rho \in [0,1)$ be such that $\chi(d_\rho) = e^{2\pi i \kappa_\rho}$. The real number $\kappa_\rho$ is called the 
{\em cusp parameter} (cf. \cite[\S 3.7]{GH}). With this setup, one has the following lemma. 
\begin{lemma}\label{fe}
Let $f\in H_k^!(\Gamma_0(N),\chi) (k \neq 1)$ and $\rho, \gamma_\rho, t_{\rho}, \kappa_\rho$ be as above. 
Then the Fourier series expansion of $f$ at the cusp $\rho$ has the following shape.
\[
(f|\gamma_{\rho})(z)=\sum_{n>>-\infty}c_f^+(n)q^{\frac{n+\kappa_{\rho}}{t_{\rho}}}+c_f^-(0)y^{1-k}q^{\frac{\kappa_{\rho}}{t_{\rho}}}+\sum_{n<<\infty}c_f^-(n)\Gamma(1-k,-4\pi ny/t_{\rho})q^{\frac{n+\kappa_{\rho}}{t_{\rho}}}.
\] 
\end{lemma}
We now briefly review the construction of Maass-Poincare series in \cite{BO}. Let $\frac{1}{2}\geq k \in \frac{1}{2}\mathbb{Z}, m \in \mathbb{N}$, and for $\chi$ a Dirichlet character modulo $N$ let $\Psi_{k,\chi}$ be the multiplier system defined in \eqref{eq 2.4}. For $L\in\mathrm{SL}(2,\mathbb{Z})$, let $\rho=\frac{a_{\rho}}{c_{\rho}}=L^{-1}(\infty)$, and $\kappa_{\rho},t_{\rho}$ be the cusp parameter and cusp width of $\rho$ with respect to $\Gamma_0(N)$. 
Let $\Gamma_{0}(N)_{\rho}$ be the stabiliser of $\rho$ in $\Gamma_0(N)$. Then Maass-Poincar\'{e} series is defined as
$$
\mathcal{F}(z,k, m,N,L,\chi):=\sum_{\gamma \in \Gamma_{0}(N)_{\rho} \backslash \Gamma_0(N)} \frac{\phi_{k}\left(\left(\frac{-m+\kappa_{\rho}}{t_{\rho}}\right)L\gamma z\right)}{j(L,\gamma z)^k j(\gamma,z)\Psi_{k,\chi}(\gamma)},
$$
where $\phi_{k}$ is given by
$$
\phi_{k}(z):=\frac{1}{\Gamma(2-k)}(-4 \pi my)^{-\frac{k}{2}} M_{\frac{k}{2}, \frac{1-k}{2}}(-4 \pi m y) e(m x).
$$
Here $M_{\mu, \nu}(w)$ are the classical $M$-Whittaker functions which are solutions of the differential equation
\[
\frac{\partial^{2} f}{\partial w^{2}}+\left(-\frac{1}{4}+\frac{\nu}{w}+\frac{\frac{1}{4}-\mu^{2}}{w^{2}}\right) f=0.
\]

By \cite[Theorem 3.2]{BO}, given a cusp $\rho$ of $\Gamma_0(N)$, we have a sequence of harmonic Maass forms $\{\varphi_{\rho,m}\}_{m=1}^{\infty}$ such that $\varphi_{\rho,m}$ decays like a cusp form at every cusp of $\Gamma_0(N)$ which is inequivalent to $\rho$ and it has a Fourier expansion of the form 
\begin{equation*}
\varphi_{\rho,m}|\gamma_{\rho}(z)=q^{\frac{-m+\kappa_{\rho}}{t_{\rho}}}+\sum_{n=0}^{\infty}c_{\varphi,m}(\rho,n)q^{\frac{n+\kappa_{\rho}}{t_{\rho}}}+\varphi_{\rho,m}^-
\end{equation*}
at the cusp $\rho.$ Here $\varphi_{\rho,m}^-$ decays like a cusp form\footnote{This follows from the fact that the shadows of Maass-Poincar\'{e} series are the classical Poincar\'{e} series which are cusp forms.} and $\gamma_{\rho}\in\mathrm{SL}(2,\mathbb{Z})$ is a uniquely chosen matrix. \par
The following Lemma is a slight modification of Lemma 2.3 of \cite{BO} and will be useful in our proof.
\begin{lemma}
Let $0> k\in\frac{1}{2}\mathbb{Z}$ and $f\in H_k^{\#}(\Gamma_0(N),\chi)$. If $\xi_k(f)\in S_{2-k}(\Gamma_0(N),\overline{\chi})$ then $f=0.$
\label{lemma 3.10}
\begin{proof}
Since $f$ has trivial principle part at every cusp of $\Gamma_0(N),$ Lemma 2.3 of \cite{BO} implies that $\xi_k(f)=0.$ This means that $f$ has trivial nonholomorphic part and hence is a holomorphic modular form of weight $0> k\in \frac{1}{2}\mathbb{Z}.$ Since there are no such non-zero forms, we conclude that $f=0.$
\end{proof}
\end{lemma}
\begin{prop}
For $0>k\in\frac{1}{2}\mathbb{Z}$, the restriction of the shadow map \[\xi_{k}:H_k^{\#}(\Gamma_0(N),\chi)\longrightarrow \frac{M_{2-k}(\Gamma_0(N),\overline{\chi})}{S_{2-k}(\Gamma_0(N),\overline{\chi})}\]
is an isomorphism of the vector spaces.
\label{thm 3.11}
\begin{proof}
We first prove surjectivity. Let \[[f]\in \frac{M_{2-k}(\Gamma_0(N),\overline{\chi})}{S_{2-k}(\Gamma_0(N),\overline{\chi})}.\] Since $\xi_k:H_k^!(\Gamma_0(N),\chi)\longrightarrow M^!_k(\Gamma_0(N),\overline{\chi})$ is surjective, let $g\in H_k^!(\Gamma_0(N),\chi)$ such that $\xi_k(g)=f.$ By \eqref{eq 2.3}, we see that the nonholomorphic part of $g$ does not have positive powers of $q$ at any cusp of $\Gamma_0(N).$ Suppose $(\rho_{i})_{i=1}^{\ell}$ be the complete set of inequivalent cusps of $\Gamma_0(N)$ with cusp width $(t_{\rho_{i}})_{i=1}^{\ell}$ and cusp parameter $(\kappa_{\rho_{i}})_{i=1}^{\ell}$ respectively. Suppose the Fourier expansion of $g$ at the cusp $\rho_i$ is given by
\[
g|\gamma_{\rho_i}(z)=\sum_{n_{\rho_i}\leq n\leq -1}c_{g}^+(\rho_i,n)q^{\frac{n+\kappa_{\rho_i}}{t_{\rho_i}}}+\sum_{n=0}^{\infty}c_{g}^+(\rho_i,n)q^{\frac{n+\kappa_{\rho_i}}{t_{\rho_i}}}+g^-(z,\rho_i)
\]
for some uniquely chosen $\gamma_{\rho_i}$.
Consider the shifted harmonic Maass form 
\[
\widehat{g}=g-\sum_{i=1}^{\ell}\widehat{\varphi}_{\rho_i}^g,\ \ {\rm where}
\]
\[
\widehat{\varphi}_{\rho_i}^g=\sum_{r=1}^{n_{\rho_i}}c_g^+(\rho_i,-n)\varphi_{\rho_i,n}.
\]
It is now easy to check that $\widehat{g}\in H_k^{\#}(\Gamma_0(N),\chi).$ Next, observe that 
\[
\xi_k(\widehat{g})=f-\xi_k\left(\sum_{i=1}^{\ell}\widehat{\varphi}_{\rho_i}^g\right)\in [f]
\]
since each $\xi_k(\varphi_{\rho_i,n})\in S_{2-k}(\Gamma_0(N),\overline{\chi}).$ Finally, an application of Lemma \ref{lemma 3.10} completes the proof of injectivity of the restriction of the shadow operator. 
\end{proof}
\end{prop}

\subsection {Proof of Theorem \ref{thm 1.5}} Since $\frac{M_{2-k}(\Gamma_0(N),\overline{\chi})}{S_{2-k}(\Gamma_0(N),\overline{\chi})}$ is isomorphic to the Eisenstein space $\mathfrak{E}_{2-k}(\Gamma_0(N),\overline{\chi})$, by using Corollary  \ref{cor 4.3} we have that $H_{k}^{\#}(\Gamma_0(N),\chi)\cong \mathfrak{E}_{k}^{\#}(\Gamma_0(N),\chi)$. Since $\mathfrak{E}_{k}^{\#}(\Gamma_0(N),\chi)\subseteq H_{k}^{\#}(\Gamma_0(N),\chi)$ we have that $\mathfrak{E}_{k}^{\#}(\Gamma_0(N),\chi)= H_{k}^{\#}(\Gamma_0(N),\chi).$ The dimension formula follows from Corollary \ref{cor 4.4}.

\end{document}